\documentclass[a4paper,12pt]{article}

\usepackage{amsthm}
\usepackage{amsmath,amssymb,latexsym,amsfonts,mathrsfs}

\renewcommand{\thefootnote}{(\arabic{footnote})}

\theoremstyle{plain}
\newtheorem{Thm}{Theorem}[section]
\newtheorem{Lem}[Thm]{Lemma}
\newtheorem{Prop}[Thm]{Proposition}
\newtheorem{Cor}[Thm]{Corollary}
\theoremstyle{definition}
\newtheorem{Def}[Thm]{Definition}
\newtheorem{Rem}[Thm]{Remark}

\newcommand{\Proof}[2][Proof]{\begin{proof}[{#1}] #2 \end{proof}}

\numberwithin{equation}{section}

\makeatletter
\renewcommand\section{\@startsection {section}{1}{\z@}%
                                   {-3.5ex \@plus -1ex \@minus -.2ex}%
                                   {2.3ex \@plus.2ex}%
                                   {\normalfont\large\bf}}
\makeatother

\makeatletter
\renewcommand\subsection{\@startsection {subsection}{1}{\z@}%
                                   {-3.5ex \@plus -1ex \@minus -.2ex}%
                                   {2.3ex \@plus.2ex}%
                                   {\normalfont\normalsize\bf}}
\makeatother

\renewcommand{\d}{{\rm d}} 

\renewcommand{\tilde}{\widetilde}

\newcommand{\bC}{\ensuremath{\mathbb{C}}}
\newcommand{\bE}{\ensuremath{\mathbb{E}}}
\newcommand{\bN}{\ensuremath{\mathbb{N}}}
\newcommand{\bP}{\ensuremath{\mathbb{P}}}
\newcommand{\bQ}{\ensuremath{\mathbb{Q}}}
\newcommand{\bR}{\ensuremath{\mathbb{R}}}
\newcommand{\bT}{\ensuremath{\mathbb{T}}}
\newcommand{\bX}{\ensuremath{\mathbb{X}}}

\newcommand{\bZ}{\ensuremath{\mathbb{Z}}}

\newcommand{\cA}{\ensuremath{\mathcal{A}}}
\newcommand{\cB}{\ensuremath{\mathcal{B}}}
\newcommand{\cC}{\ensuremath{\mathcal{C}}}
\newcommand{\cF}{\ensuremath{\mathcal{F}}}
\newcommand{\cG}{\ensuremath{\mathcal{G}}}
\newcommand{\cK}{\ensuremath{\mathcal{K}}}
\newcommand{\cN}{\ensuremath{\mathcal{N}}}
\newcommand{\cP}{\ensuremath{\mathcal{P}}}
\newcommand{\cS}{\ensuremath{\mathcal{S}}}

\newcommand{\vx}{\ensuremath{\mbox{{\boldmath $x$}}}}

\newcommand{\vw}{\ensuremath{\mbox{{\boldmath $w$}}}}

\newcommand{\rbra}[1]{\!\left( #1 \right)} 
\newcommand{\cbra}[1]{\!\left\{ #1 \right\}} 
\newcommand{\sbra}[1]{\!\left[ #1 \right]} 

\newcommand{\dist}{\stackrel{{\rm d}}{=}}
\newcommand{\cdist}{\stackrel{{\rm d}}{\longrightarrow}}

\newcommand{\tend}[2]{\mathrel{\mathop{\longrightarrow}\limits^{#1}_{#2}}}

\newcommand{\pmat}[1]{\begin{pmatrix} #1 \end{pmatrix}} 

\usepackage[usenames]{color}

\begin{document}
\begin{center}
{\Large \bf Resolution of sigma-fields 
for multiparticle finite-state action evolutions with infinite past}
\footnote{This research was supported by RIMS and by ISM.}
\end{center}

\begin{center}
Yu Ito\footnote{
Department of Mathematics, Faculty of Science, Kyoto Sangyo University, Kyoto, JAPAN.}\footnote{
The research of this author was supported by 
JSPS KAKENHI Grant Number JP18K13431.}, 
Toru Sera\footnote{
Department of Mathematics, Graduate School of Science, Osaka University, Osaka, JAPAN.
Research Fellow of Japan Society for the Promotion of Science. }\footnote{
The research of this author was supported by 
JSPS KAKENHI Grant Numbers JP19J11798 and JP21J00015.} and 
\renewcommand{\thefootnote}{(\fnsymbol{footnote})}
Kouji Yano\footnote[1]{Corresponding author. Email: {\tt kyanomath@gmail.com}}\renewcommand{\thefootnote}{(\arabic{footnote})}\footnote{
Graduate School of Science, Kyoto University, Kyoto, JAPAN.}\footnote{
The research of this author was supported by 
JSPS KAKENHI grant no.'s JP19H01791, JP19K21834 and JP18K03441 
and by JSPS Open Partnership Joint Research Projects grant no. JPJSBP120209921.} 
\end{center}

\begin{center}
{\small \today}
\end{center}

\begin{abstract}
For multiparticle finite-state action evolutions, 
we prove that the observation $ \sigma $-field 
admits a resolution involving a third noise 
which is generated by a random variable with uniform law. 
The Rees decomposition from the semigroup theory 
and the theory of infinite convolutions 
are utilized in our proofs. 
\end{abstract}

\section{Introduction}

Let us consider the stochastic recursive equation 
\begin{align}
X_k = N_k X_{k-1} 
\quad \text{for $ k \in \bZ $} , 
\label{eq: 1}
\end{align}
which we call the \emph{action evolution}, where 
the \emph{observation} $ X = (X_k)_{k \in \bZ} $ taking values 
in a measurable space $ V $ 
evolves from $ X_{k-1} $ to $ X_k $ at each time $ k $ 
being acted by a random mapping $ N_k $ of $ V $. 
Here we mean by $ N_k X_{k-1} $ the evaluation $ N_k(X_{k-1}) $ 
of a random mapping $ N_k $ at $ X_{k-1} $; 
we always abbreviate the parentheses to write $ fv $ simply for the evaluation $ f(v) $. 
As our processes are indexed by $ \bZ $, 
the state $ X_k $ we observe at time $ k $ is a result 
after a long time has passed.

We would like to clarify the structure of 
the full noise $ \cF^{X,N}_k = \sigma(X_j,N_j : j \le k) $ 
and the observation noise $ \cF^X_k = \sigma(X_j : j \le k) $. 
For families of events, 
we write $ \cA \vee \cB := \sigma(\cA \bigcup \cB) $. 
For $ \sigma $-fields, 
we say $ \cF \subset \cG $ a.s. (resp. $ \cF = \cG $ a.s.) 
if $ \cF \subset \cG \vee \cN $ (resp. $ \cF \vee \cN = \cG \vee \cN $) 
with $ \cN $ being the family of null events. 
By iterating the equation \eqref{eq: 1}, 
we have $ X_k = N_k N_{k-1} \cdots N_{j+1} X_j $ a.s. for $ j<k $. 
One may then expect that, for any $ k \in \bZ $, 
\begin{align}
\cF^{X,N}_k = \bigcap_{j<k} \rbra{ \cF^N_k \vee \cF^X_j } 
\stackrel{?}{=} \cF^N_k \vee \rbra{ \bigcap_{j<k} \cF^X_j } 
= \cF^N_k \vee \cF^X_{-\infty } 
\quad \text{a.s.}, 
\label{eq: intro inclusion}
\end{align}
and may conclude that the full noise $ \cF^{X,N}_k $ 
can be known by the \emph{driving noise} $ \cF^N_k := \sigma(N_j : j \le k) $ 
together with the \emph{remote past noise} $ \cF^X_{-\infty } := \bigcap_k \cF^X_k $, 
which plays the role of the initial noise at time $ -\infty $. 
But the a.s. identity $ \stackrel{?}{=} $ in \eqref{eq: intro inclusion} 
is false in general; see \cite[(1) of Remark 1.4]{MR3009741} 
for erroneous discussions by Kolmogorov and Wiener. 
We must refer to \cite[Section 2.5]{CY} for careful treatment 
of exchanging the order of supremum and intersection between $ \sigma $-fields.

\subsection{Action evolutions and resolution of the observation}

We would like to reveal the hidden extra noise. 
To this end let us introduce some terminology. 
The action evolution proposed in \eqref{eq: 1} is formulated as follows. 

\begin{Def}
Let $ \mu $ be a probability on a measurable space $ V^V $ of mappings of $ V $ into itself 
and call it a \emph{mapping law} on $ V^V $. 
A (mono-particle) \emph{$ \mu $-evolution} is a pair $ (X,N) $ of 
a $ V $-valued process $ X = (X_k)_{k \in \bZ} $ 
and an iid $ V^V $-valued process $ N = (N_k)_{k \in \bZ} $ 
defined on a probability space $ (\Omega,\cF,\bP) $ such that 
the following hold for each $ k \in \bZ $: 
\\ \quad (i) 
$ X_k = N_k X_{k-1} $ holds a.s.; 
\\ \quad (ii) 
$ N_k $ is independent of $ \cF^{X,N}_{k-1} := \sigma(X_j,N_j : j \le k-1) $; 
\\ \quad (iii) 
$ N_k $ has law $ \mu $. 
\end{Def}

It is easy to see that $ (X,N) $ is a $ \mu $-evolution if and only if 
the Markov property 
\begin{align}
\bP \rbra{ (X_k,N_k) \in B \mid \cF^{X,N}_{k-1} } 
= Q_\mu \rbra{ X_{k-1}; B \Big. } 
, \quad k \in \bZ , \ B \subset V \times V^V 
\label{eq: MP}
\end{align}
holds with the joint transition probability: 
\begin{align}
Q_\mu \rbra{ x; B \Big. } = \mu \cbra{ f : (fx,f) \in B \Big. } 
, \quad x \in V , \ B \subset V \times V^V , 
\label{}
\end{align}
i.e., $ Q_\mu(x;\cdot) $ is the image measure of $ \mu $ by the map $ f \mapsto (fx,f) $. 
If $ (X,N) $ is a $ \mu $-evolution, then 
the marginal process $ X $ satisfies the Markov property 
\begin{align}
\bP \rbra{ X_k \in A \mid \cF^X_{k-1} } 
= P_\mu(X_{k-1};A) 
, \quad k \in \bZ , \ A \subset V 
\label{}
\end{align}
with the marginal transition probability: 
\begin{align}
P_\mu(x;A) = \mu \cbra{ f : fx \in A } 
, \quad A \subset V , 
\label{eq: mtp}
\end{align}
i.e., $ P_\mu(x;\cdot) $ is the image measure of $ \mu $ by the map $ f \mapsto fx $. 

It is easy to see by definition that, if two $ \mu $-evolutions $ (X,N) $ and $ (X',N') $ satisfy 
$ X_k \dist X'_k $ for $ k \in \bZ $, then $ (X,N) \dist (X',N') $. 
This shows that the joint law of the $ \mu $-evolution $ (X,N) $ is determined 
by the family of marginal laws $ \{ \bP(X_k \in \cdot) : k \in \bZ \} $. 

Unfortunately, it seems difficult to develop 
a complete investigation of the structure of the observation noise 
for a mono-particle $ \mu $-evolution. 
So we introduce a multi-particle counterpart. 
For a mapping $ f:V \to V $ and a vector $ \vx = (x^1,\ldots,x^m) \in V^m $, 
we understand that $ f $ operates $ \vx $ componentwise, i.e., $ f \vx = (fx^1,\ldots,fx^m) $. 

\begin{Def} \label{def: mu-evol}
Let $ \mu $ be a mapping law on $ V^V $. 
An \emph{$ m $-particle $ \mu $-evolution} is a $ \mu $-evolution $ (\bX,N) $ 
with $ \bX = (\bX_k)_{k \in \bZ} $ taking values in $ V^m $; 
precisely, the following hold for each $ k \in \bZ $: 
\\ \quad (i) 
$ \bX_k = N_k \bX_{k-1} $ holds a.s., i.e., 
$ X^i_k = N_k X^i_{k-1} $ holds a.s. for $ i=1,\ldots,m $; 
\\ \quad (ii) 
$ N_k $ is independent of $ \cF^{\bX,N}_{k-1} := \sigma(\bX_j,N_j : j \le k-1) $; 
\\ \quad (iii) 
$ N_k $ has law $ \mu $. 
\end{Def}

We will see in Proposition \ref{prop: num of distinct} that 
the number of distinct states among $ \{ X^1_k,\ldots,X^m_k \} $ 
does not depend upon $ k \in \bZ $ a.s.

We now propose our resolution problems. 

\begin{Def} \label{def: rfn}
For an $ m $-particle $ \mu $-evolution $ (\bX,N) $, 
a \emph{third noise} is a sequence of random variables $ (U_k)_{k \in \bZ} $ 
such that the following hold for each $ k \in \bZ $: 
\\ \quad (i) 
the identity 
$ \cF^{\bX,N}_k = \cF^N_k \vee \cF^{\bX}_{-\infty } \vee \sigma(U_k) $ holds a.s.; 
\\ \quad (ii) 
the three $ \sigma $-fields $ \cF^N_k $, $ \cF^{\bX}_{-\infty } $ and $ \sigma(U_k) $ 
are independent. 
\\
The identity in Condition (i) will be called the \emph{resolution of the full noise}. 
\end{Def}

We remark that the third noise is not an innovation. 
Note that, if $ (U_k)_{k \in \bZ} $ is a third noise, 
then, for any $ k_0 \in \bZ $, the stopped sequence $ (U_{k \wedge k_0})_{k \in \bZ} $ 
is also a third noise, since $ X_k = N_k N_{k-1} \cdots N_{k_0+1} X_{k_0} $ for $ k>k_0 $. 
We may suggest that 
the third noise emerges in the remote past and does not increase as time passes. 

Provided that we can find a random variable $ \xi $ such that 
$ \cF^{\bX}_{-\infty } = \sigma(\xi) $, then, for every $ k \in \bZ $, 
the identity $ \cF^{\bX,N}_k = \cF^N_k \vee \cF^{\bX}_{-\infty } \vee \sigma(U_k) $ holds 
if and only if 
there exist measurable mappings $ F $ and $ G $ such that 
$ (\bX_k,\bX_{k-1},\ldots) = F(N_k,N_{k-1},\ldots, \xi, U_k) $ 
and $ U_k = G(\bX_k,\bX_{k-1},\ldots,N_k,N_{k-1},\ldots) $.

We also propose a finer resolution. 

\begin{Def} \label{def: ro}
For an $ m $-particle $ \mu $-evolution, 
a \emph{reduced driving noise} is 
a sequence of $ \sigma $-fields $ (\cG^N_k)_{k \in \bZ} $ 
accompanying with a sequence of random variables $ (U_k)_{k \in \bZ} $ 
such that the following hold for each $ k \in \bZ $: 
\\ \quad (i) 
the identity 
$ \cF^{\bX}_k = \cG^N_k \vee \cF^{\bX}_{-\infty } \vee \sigma(U_k) $ holds a.s.; 
\\ \quad (ii) 
$ \cG^N_k \subset \cF^N_k $ holds a.s.; 
\\ \quad (iii) 
the three $ \sigma $-fields $ \cF^N_k $, $ \cF^{\bX}_{-\infty } $ and $ \sigma(U_k) $ 
are independent. 
\\
The identity in Condition (i) will be called the \emph{resolution of the observation}. 
\end{Def}

For every $ k \in \bZ $, provided that we can find random variables $ \xi $ 
and $ \zeta_k,\zeta_{k-1},\ldots $ such that 
$ \cF^{\bX}_{-\infty } = \sigma(\xi) $ and 
$ \cG^N_k = \sigma(\zeta_k,\zeta_{k-1},\ldots) $, then 
the identity $ \cF^{\bX}_k = \cG^N_k \vee \cF^{\bX}_{-\infty } \vee \sigma(U_k) $ holds 
if and only if there exist measurable mappings $ F $ and $ G $ such that 
$ (\bX_k,\bX_{k-1},\ldots) = F(\zeta_k,\zeta_{k-1},\ldots, \xi, U_k) $ 
and $ U_k = G(\bX_k,\bX_{k-1},\ldots) $. 

We remark that, 
if we have the resolution of the observation in the sense of Definition \ref{def: ro}, 
then we also have the resolution of the full noise in the sense of Definition \ref{def: rfn}, 
and the sequence $ (U_k)_{k \in \bZ} $ in Definition \ref{def: ro} 
becomes a third noise in the sense of Definition \ref{def: rfn}.

In this paper, we shall give a general result of resolution of the observation 
for multiparticle action evolutions when the state space $ V $ is a finite set.

\subsection{Infinite convolutions on finite semigroups}

For our purpose we need several known facts from the theory of semigroups, 
which we recall without proofs. 
We may consult a celebrated textbook \cite{HM} 
and a concise expository note \cite{Yinf} for the details. 

In what follows we assume $ S $ be a finite semigroup 
and we denote the set of all idempotents in $ S $ 
by $ E(S) = \{ f \in S : f^2 = f \} $. 
For $ A,B \subset S $ and $ f \in S $, we write 
$ AB = \{ ab : a \in A , \ b \in B \} $ and $ Af = \{ af : a \in A \} $, etc. 
We say that $ S $ is \emph{completely simple} if $ S $ has no proper ideal, 
i.e., $ \emptyset \neq IS \cup SI \subset I \subset S $ implies $ I=S $, 
and if there exists $ e \in E(S) $ which is \emph{primitive}, 
i.e., $ ef = fe = f \in E(S) $ implies $ f=e $.

\begin{Prop}[Rees decomposition] \label{prop: RD}
Suppose $ S $ be a completely simple finite semigroup 
with a primitive idempotent $ e $. 
Set $ L = E(Se) $, $ G = eSe $ and $ R = E(eS) $. 
Then the following hold: 
\begin{enumerate}
\item 
$ G $ is a group whose identity is $ e $. 
\item 
$ eL = Re = \{ e \} $. 
\item 
$ S = LGR $. 
\item 
The product mapping 
$ \psi: L \times G \times R \ni (f,g,h) \mapsto fgh \in S $ is bijective 
and its inverse is given as 
\begin{align}
\psi^{-1}(z) \rbra{ =: (z^L,z^G,z^R) } = (ze (eze)^{-1}, eze, (eze)^{-1} ez) . 
\label{}
\end{align}
\end{enumerate}
\end{Prop}

The proof of Proposition \ref{prop: RD} 
can be found, e.g., in \cite[Theorem 1.1]{HM}. 
The product decomposition $ S = LGR $ will be called 
the \emph{Rees decomposition} of $ S $ at $ e $, 
and $ G $ will be called the \emph{group factor}. 
Note that the Rees decomposition depends upon the choice of a primitive idempotent; 
this is why we clarify the choice of the primitive idempotent 
by saying ``at $ e $" whenever we call the Rees decomposition. 

Note by definition that $ RL \subset eSSe \subset eSe = G $ 
and by the product bijectivity that 
$ \psi^{-1}((fgh)(f'g'h')) = (f,ghf'g',h') $. 
It is obvious that the product $ z = fgh \in S $ is idempotent if and only if $ hf = g^{-1} $. 
We notice that \emph{all idempotents of $ S $ are primitive}; 
in fact, if $ e' = f'g'h' \in E(S) $ and $ z = fgh \in S $ 
satisfies $ e'z = ze' = z \in E(S) $, 
then we have $ f'=f $ and $ h'=h $ by the product bijectivity 
and thus we have $ g' = (h'f')^{-1} = (hf)^{-1} = g $, 
which shows $ e' = z $ so that $ e' $ is also primitive.

A subset $ K $ of $ S $ is called a \emph{kernel} of $ S $ 
if it is a minimal ideal of $ S $, i.e., 
$ K $ is an ideal of $ S $ and contains no proper ideal of $ S $. 
Note that a kernel $ K $ of $ S $ is then automatically the least ideal of $ S $; 
indeed, for any ideal $ I $ of $ S $, we see that $ K \cap I $ 
is an ideal of $ S $ containing $ KI $ 
and is contained in $ K $, which implies $ K \subset K \cap I \subset I $.

\begin{Prop} \label{prop: kernel}
A finite semigroup $ S $ always contains a unique kernel $ K $. 
In addition, the kernel $ K $ is completely simple. 
Fix $ e \in E(K) $ 
and let $ K = LGR $ denote the Rees decomposition of $ K $ at $ e $. 
Then 
\begin{align}
Se = Ke = LG 
, \quad 
eSe = eKe = G 
, \quad 
eS = eK = GR . 
\label{eq: Se eSe eS}
\end{align}
\end{Prop}

\Proof{
The proof of unique existence and complete simplicity of the kernel 
can be found, e.g., in \cite[Proposition 1.7]{HM}. 
Let us prove \eqref{eq: Se eSe eS}. Since 
\begin{align}
Se = See \subset SKe \subset Ke = LGRe = LGe \subset Se , 
\label{}
\end{align}
we obtain $ Se = Ke = LGe = LG $. By a similar argument we have $ eS = GR $. 
Thus we obtain $ eSe = eLG = G $. 
}

Proposition \ref{prop: RD} is fundamental in the theory of infinite convolutions. 
Let $ \cP(S) $ denote the set of probability measures on a finite semigroup $ S $ 
and write $ \mu * \nu $ for the convolution of $ \mu $ and $ \nu $ in $ \cP(S) $: 
\begin{align}
(\mu * \nu)(A) = \sum_{f,g \in S} 1_A(fg) \mu\{ f \} \nu\{ g \} 
, \quad A \subset S . 
\label{}
\end{align}
We write $ \mu^n $ for the $ n $-fold convolution of $ \mu $: 
$ \mu^1 = \mu $ and $ \mu^n = \mu * \mu^{n-1} $ for $ n=2,3,\ldots $ 
We write $ \cS(\mu) = \{ f \in S : \mu \{ f \} > 0 \} $ for the support of $ \mu $. 
It is easy to see that $ \cS(\mu * \nu) = \cS(\mu) \cS(\nu) $ for $ \mu,\nu \in \cP(S) $. 
We write $ \omega_A $ for the uniform distribution on a finite set $ A $; 
if $ G $ is a finite group, then $ \omega_G $ is the normalized Haar measure of $ G $.

\begin{Prop}[Convolution idempotents] \label{prop: CI}
Suppose that $ \nu^2 = \nu \in \cP(S) $. 
Then $ \cS(\nu) $ is a completely simple subsemigroup of $ S $. 
Fix $ e \in E(\cS(\nu)) $ and 
take $ L = E(\cS(\nu) e) $, $ G = e \cS(\nu) e $ and $ R = E(e \cS(\nu)) $ 
so that $ \cS(\nu) = LGR $ gives the Rees decomposition of $ \cS(\nu) $ at $ e $. 
Write $ \nu^L(\cdot) = \nu \{ z \in \cS(\nu) : z^L \in \cdot \} $ 
and $ \nu^R(\cdot) = \nu \{ z : z^R \in \cdot \} $, 
i.e., $ \nu^L $ and $ \nu^R $ are the image measures of $ \nu $ 
by the maps $ z \mapsto z^L $ and $ z \mapsto z^R $, respectively. 
Then $ \nu $ has a factorization 
\begin{align}
\nu = \nu^L * \omega_G * \nu^R . 
\label{}
\end{align}
Consequently, 
if $ Z $ is a random variable whose law is $ \nu $, 
then the projections $ Z^L $, $ Z^G $ and $ Z^R $ are independent 
and $ Z^G $ is uniform on $ G $. 
\end{Prop}

The proof of Proposition \ref{prop: CI} 
can be found, e.g., in \cite[Theorem 2.2]{HM}.

The following proposition plays a key role in our analysis. 

\begin{Prop}[Infinite convolutions] \label{prop: IC}
Let $ \mu \in \cP(S) $ 
and suppose that $ S $ coincide with $ \bigcup_{n=1}^{\infty } \cS(\mu)^n $, 
the semigroup generated by $ \cS(\mu) $. 
Then the following hold: 
\begin{enumerate}
\item 
The set of subsequential limits of $ \{ \mu^n \} $ 
is a finite cyclic group of the form 
\begin{align}
\cK := \{ \eta,\mu * \eta,\ldots,\mu^{p-1} * \eta \} 
\label{}
\end{align}
for some $ p \in \bN $, where $ \eta $ is the identity of $ \cK $ (so that $ \eta^2 = \eta $), 
$ \mu^p * \eta = \eta $ and 
$ \eta,\mu * \eta,\ldots,\mu^{p-1} * \eta $ are all different 
(and consequently $ \cK $ is a cyclic group with order $ p $). 
The support $ \cS(\eta) $ is a completely simple subsemigroup of $ S $ 
(but not in general an ideal of $ S $.) 
\item 
It holds that 
\begin{align}
\frac{1}{n} \sum_{k=1}^n \mu^k 
\tend{}{n \to \infty } \nu := \frac{1}{p} \sum_{k=0}^{p-1} \mu^k * \eta , 
\label{eq: Cesaro limit}
\end{align}
so that $ \nu^2 = \nu $. The support $ \cS(\nu) $ is the kernel of $ S $. 
\item 
Let $ e \in E(\cS(\eta)) $ be fixed. Then the Rees decompositions at $ e $ 
of $ \cS(\nu) $ and of $ \cS(\eta) $ are given as 
\begin{align}
\cS(\nu) = LGR 
\quad \quad \text{and} \quad 
\cS(\eta) = LHR , 
\label{}
\end{align}
respectively, where $ L = E(\cS(\eta)e) $, $ R = E(e\cS(\eta)) $, 
$ G = e\cS(\nu)e $ and $ H = e\cS(\eta)e $. 
Moreover, the group factor $ H $ of $ \cS(\eta) $ is a normal subgroup 
of the group factor $ G $ of $ \cS(\nu) $, and 
the convolution factorizations of $ \nu $ and $ \eta $ are given as 
\begin{align}
\nu = \eta^L * \omega_G * \eta^R 
\quad \quad \text{and} \quad 
\eta = \eta^L * \omega_H * \eta^R , 
\label{}
\end{align}
respectively, where 
$ \eta^L(\cdot) = \eta \{ z : z^L \in \cdot \} $ 
and $ \eta^R(\cdot) = \eta \{ z : z^R \in \cdot \} $, 
i.e., $ \eta^L $ and $ \eta^R $ are the image measures of $ \eta $ 
by the maps $ z \mapsto z^L $ and $ z \mapsto z^R $, respectively. 
\item 
There exists $ \gamma \in G $ such that 
$ G/H = \{ H,\gamma H,\ldots,\gamma^{p-1}H \} $ with $ \gamma^p \in H $ 
and with $ H,\gamma H,\ldots,\gamma^{p-1}H $ are all different 
(and consequently $ G/H $ is a cyclic group with order $ p $ 
and generated by $ \gamma H $). Moreover, 
\begin{align}
\mu^r * \eta = \eta^L * \delta_{\gamma^r} * \omega_H * \eta^R 
, \quad r=0,1,\ldots,p-1 , 
\label{}
\end{align}
where $ \delta_a $ stands for the Dirac mass at $ a $. 
\end{enumerate}
\end{Prop}

The proof of Proposition \ref{prop: IC} 
can be found in \cite{Yinf}; see also \cite[Theorem 2.7]{HM}.

\begin{Rem} \label{rem: mupn to eta}
By (i) and (iv) of Proposition \ref{prop: IC}, we see that 
$ \mu^{pn} \to \eta $, $ \mu^{pn+1} \to \mu * \eta, \ldots $, 
$ \mu^{pn+p-1} \to \mu^{p-1} * \eta $. 
To see this fact, it suffices to prove $ \mu^{pn} \to \eta $. 
Since $ \cP(S) $ is compact, it suffices to prove that 
$ \eta $ is the unique cluster point of $ \{ \mu^{pn} \} $. 
Suppose that a subsequence $ \mu^{pn(k)} $ converge 
to an element of $ \cK $, say $ \mu^r * \eta $ for $ r=0,1,\ldots,p-1 $. 
Since $ \mu^{pn(k)} * \eta = \eta $, we have $ \mu^r * \eta * \eta = \eta $, 
which implies $ r=0 $. 
We thus obtain $ \mu^{pn(k)} \to \eta $. 
\end{Rem}

\begin{Rem}
It may be useful to notice some connection between Proposition \ref{prop: IC} 
and random walks generated by $ \mu $. 
By \cite[Proposition 3.6]{HM} and by \cite[Corollary 3.1, Theorem 3.2 and Proposition 3.4]{HM}, 
we have the following facts: 
\\ $ \bullet $ 
The set $ \cS(\nu) $, which is the kernel of $ S $, is also the set of all recurrent 
elements for the unilateral, the bilateral or the mixed random walk generated by $ \mu $. 
\\ $ \bullet $ 
For any $ z \in \cS(\nu) $, the greatest common divisor of 
$ \{ n-m : \mu^n \{ z \} > 0 , \ \mu^m \{ z \} > 0 \} $ 
coincides with $ p $. 
\end{Rem}

\subsection{The semigroup consisting of mappings}

Let $ V $ be a non-empty finite set 
and let $ V^V $ denote the set of mappings of $ V $ into itself. 
Note that $ V^V $ is also a finite semigroup with respect to composition 
as its product structure. 
In this concrete settings 
we recall several facts 
about the description of the kernel and the Rees decomposition 
(see e.g. \cite[Example 1.1 and Proposition 1.8]{HM}). 
For $ f \in V^V $, we write $ \pi(f) := \{ f^{-1} \{ v \} : v \in V \} $ 
for the partition of $ S $ generated by the preimages of $ f $. 

\begin{Prop} \label{prop: mapping semigroup} 
Let $ S $ be a subsemigroup of $ V^V $ and denote 
\begin{align}
m_S = \min \{ \# (fV) : f \in S \} , 
\label{}
\end{align}
where $ \#(A) $ denotes the number of elements of a set $ A $. 
\begin{enumerate}

\item 
The kernel $ K $ of $ S $ is the set of all mappings in $ S $ with minimal rank, namely: 
\begin{align}
K = \{ f \in S : \# (fV) = m_S \} . 
\label{eq: def of K}
\end{align}

\item 
$ e $ is an idempotent if and only if it is identity on $ eV $. 

\item 
$ e $ is a primitive idempotent of $ S $ if and only if $ e \in E(K) $.

\end{enumerate}
For a fixed $ e \in E(K) $, 
the Rees decomposition $ K=LGR $ at $ e $ may be characterized as follows: 
\begin{enumerate}
\addtocounter{enumi}{3}

\item 
$ Se $ is the set of all $ f $ in $ S $ such that $ \pi(f) = \pi(e) $. 
Consequently, 
$ L = E(Se) $ is the set of all idempotents $ f $ in $ S $ such that $ \pi(f) = \pi(e) $. 

\item 
$ eS $ is the set of all $ f $ in $ S $ such that $ fV = eV $. 
Consequently, 
$ R = E(eS) $ is the set of all idempotents $ f $ in $ S $ such that $ fV = eV $. 

\item 
$ G $ is the set of all $ f \in S $ such that 
$ \pi(f) = \pi(e) $ and $ fV = eV $.
\end{enumerate}
\end{Prop}

For convenience of the readers, 
the proof of Proposition \ref{prop: mapping semigroup} will be given in the Appendix.

\subsection{Main result} \label{subsec: main results}

Let $ V $ be a non-empty finite set 
and let $ V^V $ denote the set of mappings of $ V $ into itself. 
For $ f \in V^V $ and $ \vx = (x^1,\ldots,x^m) \in V^m $, 
we understand $ f \vx = (fx^1,\ldots,fx^m) $. 
For $ \mu \in \cP(V^V) $ and $ \Lambda \in \cP(V^m) $, 
we define $ \mu * \Lambda \in \cP(V^m) $ as 
\begin{align}
(\mu * \Lambda)(A) = \sum_{f \in V^V} \sum_{\vx \in V^m} 1_A(f \vx) \mu \{ f \} \Lambda \{ \vx \} 
, \quad A \subset V^m . 
\label{}
\end{align}
Note that, for independent random variables $ F $ and $ \bX $ 
whose laws are $ \mu $ and $ \Lambda $, respectively, 
the law of $ F \bX $ is $ \mu * \Lambda $. 
Denote 
\begin{align}
V^m_\times = \{ \vx = (x^1,\ldots,x^m) \in V^m : \text{$ x^1,\ldots,x^m $ are distinct} \} . 
\label{}
\end{align}

\begin{Prop} \label{prop: contr W}
Let $ \mu \in \cP(V^V) $ 
and set $ S = \bigcup_{n=1}^{\infty } \cS(\mu)^n $, 
the semigroup generated by $ \cS(\mu) $. 
We apply Proposition \ref{prop: IC} and adopt its notations. 
Denote 
\begin{align}
m_\mu = m_S = \min \{ \# (fV) : f \in S \} 
\label{}
\end{align}
and define 
\begin{align}
W_\mu = \{ \vx \in V^{m_\mu}_\times : 
\text{$ f\vx \in V^{m_\mu}_\times $ for all $ f \in S $} \} . 
\label{eq: Wmu}
\end{align}
(Note that $ W_\mu $ is not empty; in fact, if we write $ eV = \{ x_1,\ldots,x_{m_\mu} \} $, 
then $ (x_1,\ldots,x_{m_\mu}) \in W_\mu $.) 
Let $ W $ be an arbitrary minimal subset of $ eW_\mu $ such that $ eW_\mu = GW $. 
Then the following hold: 
\begin{enumerate}

\item 
$ W_\mu = LGW $. 

\item 
The product mapping $ L \times G \times W \ni (f,g,\vw) \mapsto fg\vw \in W_\mu $ 
is bijective. 
Its inverse will be denoted by 
$ \vx \mapsto (\vx^L,\vx^G,\vx^W) $. 

\item 
Let $ \Lambda \in \cP(V^{m_\mu}_\times) $. 
Then $ \Lambda $ is \emph{$ \mu $-invariant}, i.e., $ \Lambda = \mu * \Lambda $, 
if and only if 
$ \Lambda = \eta^L * \omega_G * \Lambda_W $ 
for some $ \Lambda_W \in \cP(W) $. 

\end{enumerate}
\end{Prop}

The proof of Proposition \ref{prop: contr W} will be given in Section \ref{sec: F-cliques}.

If an $ m $-particle $ \mu $-evolution $ (\bX,N) $ is \emph{stationary}, 
i.e., $ (\bX_{\cdot + 1},N_{\cdot + 1}) \dist (\bX,N) $, 
then the sequence $ \bX $ has a common law which is $ \mu $-invariant. 
Conversely, if $ \Lambda \in \cP(V^m) $ is $ \mu $-invariant, 
then there exists a stationary $ m $-particle $ \mu $-evolution $ (\bX,N) $ 
such that the sequence $ \bX $ has $ \Lambda $ as its common law. 

In Proposition \ref{prop: IC}, 
we write $ C = \{ e,\gamma,\ldots,\gamma^{p-1} \} $ so that $ CH = \bigcup_{c \in C} cH = G $. 
By definition of $ C $ and $ H $, 
we see that the product mapping 
$ C \times H \ni (c,h) \mapsto ch \in G $ is bijective. 
Its inverse will be denoted by $ g \mapsto (g^C,g^H) $. 

We now state our main theorem, which will be proved in Section \ref{sec: proof}. 

\begin{Thm} \label{thm: main}
Suppose the same assumptions of Proposition \ref{prop: contr W} be satisfied. 
Suppose that $ \Lambda \in \cP(V^{m_\mu}_\times) $ be $ \mu $-invariant 
and let $ (\bX,N) $ be a stationary $ m_\mu $-particle $ \mu $-evolution 
such that the sequence $ \bX $ has $ \Lambda $ as its common law. 
Then the following hold: 
\begin{enumerate}

\item 
For any fixed $ k \in \bZ $, 
it holds that $ \bX_k \in W_\mu = LGW $ a.s., 
$ \bX_k^L \dist \eta^L $, $ \bX_k^G \dist \omega_G $, $ \bX_k^W \dist \Lambda_W $ 
and the three random variables $ \bX_k^L $, $ \bX_k^G $ and $ \bX_k^W $ 
are independent. 
(Note that the two processes $ \bX^L $ and $ \bX^G $ are \emph{not} independent 
in general; see \eqref{eq: not indep}.) 

\item 
$ \bX_k^G = (\gamma^k Y_C)^C U_k $ a.s. for $ k \in \bZ $ 
for some $ C $-valued random variable $ Y_C $ and some $ H $-valued random variables $ U_k $ 
such that $ U_k $ is uniform on $ H $. 

\item 
$ \bX_k^W = \bZ_W $ a.s. for $ k \in \bZ $ 
for some $ W $-valued random variable $ \bZ_W $. 

\item 
If we write $ M^G_j := \bX^G_j (\bX^G_{j-1})^{-1} $ for $ j \in \bZ $ 
and $ M^G_{k,j} := \bX^G_k (\bX^G_j)^{-1} = M^G_k M^G_{k-1} \cdots M^G_{j+1} $ for $ j \le k $, 
we have the following factorization: 
\begin{align}
\bX_j = \bX_j^L (M^G_{k,j})^{-1} (\gamma^k Y_C)^C U_k \bZ_W 
\quad \text{a.s. for $ j \le k $}. 
\label{eq: bXj factorization}
\end{align}

\item 
A resolution of the observation holds in the sense that 
\begin{align}
\cF^{\bX}_k = \cG^N_k \vee \cF^{\bX}_{-\infty } \vee \sigma(U_k) 
\quad \text{a.s.}, 
\label{eq: resol}
\end{align}
where 
\begin{align}
\cG^N_k = \sigma \rbra{ \bX^L_j , \ M^G_j : j \le k } \subset \cF^N_k (\subset \sigma(N)) 
\quad \text{a.s.} , 
\label{eq: cGNk cFNk}
\end{align}
\begin{align}
\text{the three $ \sigma $-fields $ \sigma(N) $, 
$ \cF^{\bX}_{-\infty } $ and $ \sigma(U_k) $ are independent} 
\label{eq: main indep}
\end{align}
and 
\begin{align}
\cF^{\bX}_{-\infty } 
= \sigma(Y_C , \ \bZ_W) 
\quad \text{a.s.} 
\label{eq: bX-infty}
\end{align}

\item 
$ Y_C \dist \omega_C $ 
and $ \bZ_W \dist \Lambda_W $, 
where $ \omega_C $ denotes the uniform distribution on the set $ C $. 
It holds that $ Y_C $ and $ \bZ_W $ are independent. 

\end{enumerate}
\end{Thm}

We shall show in Section \ref{sec: non-stationary} 
that the non-stationary case can be reduced to the stationary case 
and satisfies Properties (i)-(v) of Theorem \ref{thm: main}. 

Note that, if we represent $ Y_C = \gamma^{R_C} $ with $ R_C \in \{ 0,1,\ldots,p-1 \} $, then 
\begin{align}
(\gamma^k Y_C)^C 
= (\gamma^{k+R_C})^C = \gamma^{(k+R_C) \mathop{\rm mod} p} 
\quad \text{for $ k \in \bZ $}, 
\label{}
\end{align}
since $ \gamma^p \in H $.

The following corollary to Proposition \ref{prop: contr W} 
ensures that a stationary mono-particle $ \mu $-evolution 
can always be extended to a stationary $ m_\mu $-particle $ \mu $-evolution. 

\begin{Cor} \label{cor: extension}
Let $ \lambda \in \cP(V) $ and suppose that $ \lambda $ be $ \mu $-invariant, 
i.e., $ \lambda = \mu * \lambda $. 
Then there exists $ \Lambda \in \cP(V^{m_\mu}_\times) $ such that 
$ \Lambda $ is $ \mu $-invariant and 
its marginal in the first coordinate equals to $ \lambda $, i.e., 
\begin{align}
\Lambda \{ (x^1,\ldots,x^{m_\mu}) \in V^{m_\mu}_\times : x^1 = v \} = \lambda \{ v \} 
, \quad v \in V . 
\label{}
\end{align}
\end{Cor}

The proof of Corollary \ref{cor: extension} 
will be given in Section \ref{sec: F-cliques}. 
Unfortunately, the resolution of the observation obtained in Theorem \ref{thm: main} 
for an $ m_\mu $-paricle $ \mu $-evolution 
does not imply that for a mono-paricle $ \mu $-evolution, 
which will be illustrated in Section \ref{sec: ex}.

\subsection{Historical remarks} \label{sec: HR}

The theories of Rees decomposition, convolution idempotents 
and infinite convolutions for finite semigroups are very old results 
and have nowadays been generalized to topological semigroups; 
see the textbook \cite[Chapters 1 and 2]{HM} for the details. 
In particular, Proposition \ref{prop: IC}, which plays a fundamental tool for our results, 
dates back to Rosenblatt \cite{MR0118773}, Collins \cite{MR137789} and Schwarz \cite{MR169969}.

Inspired by Tsirelson \cite{MR0375461} of a stochastic differential equation, 
Yor \cite{MR1147613} has made a thorough study of the action evolution 
$ X_k = N_k X_{k-1} $ when both $ X $ and $ N $ take values in the torus 
$ \bT = \{ z \in \bC : |z|=1 \} $ and $ N $ is not necessarily iid, 
where we understand $ N_k X_{k-1} $ as the usual product 
between two complex values. 
He obtained a general result of the resolution of the observation. 
Hirayama and Yano \cite{MR2653259} generalized Yor's results 
for the state space being a compact group. 
In these results the third noise is generated by a random variable 
with uniform law on a subgroup of the state space group. 
See also \cite{MR3374628} for a survey of this topic.

If $ m_\mu = 1 $, then it is obvious that 
a stationary mono-particle $ \mu $-evolution $ (X,N) $ satisfies 
$ \cF^X_k \subset \cF^N_k $ a.s. for all $ k $. 
Yano \cite{MR3023844} proved its converse: 
if a stationary mono-particle $ \mu $-evolution $ (X,N) $ satisfies 
$ \cF^X_k \subset \cF^N_k $ a.s. for all $ k $, 
then $ m_\mu = 1 $. For the proof, 
he proved existence of a non-trivial third noise when $ m_\mu \ge 2 $. 
He utilized several notions from the \emph{road coloring theory}; 
for the details see Trahtman \cite{MR2534238} and the references therein.

Brossard--Leuridan \cite{MR2308593} have studied 
Markov chains indexed by $ \bZ $, which can be regarded 
as $ \mu $-evolutions on general state spaces. 
Let us pick up \cite[Theorem 3]{MR2308593}, 
whose conclusion applied to our $ \mu $-evolutions is as follows: 
\emph{For a stationary $ m_\mu $-particle $ \mu $-evolution $ (\bX,N) $ 
such that $ \bX $ has $ \Lambda $ as its common law, 
the conditional law of $ \bX $ given $ \sigma(N) \vee \cF^{\bX}_{-\infty } $ 
is uniform on some random finite set.} 
Let us characterize the random finite set. 
By Theorem \ref{thm: main}, we see that 
the conditional law $ \bP(\bX \in \cdot \mid \sigma(N) \vee \cF^{\bX}_{-\infty }) $ 
is uniform on the random finite set $ \{ \bX^h : h \in H \} \subset (V^{m_\mu})^{\bZ} $, 
where the process $ \bX^h = (\bX^h_j)_{j \in \bZ} $ is defined as 
\begin{align}
\bX^h_j = 
\begin{cases}
\bX^L_j (M^G_{0,j})^{-1} Y_C h \bZ_W & (j \le 0) \\
N_j N_{j-1} \cdots N_1 \bX^L_0 Y_C h \bZ_W & (j \ge 1) . 
\end{cases}
\label{}
\end{align}
In particular, the conditional law 
$ \bP(\bX_0 \in \cdot \mid \sigma(N) \vee \cF^{\bX}_{-\infty }) $ 
is uniform on the random finite set 
$ \{ \bX^L_0 Y_C h \bZ_W : h \in H \} \subset V^{m_\mu} $. 
In Remark \ref{rem: BL}, 
we will discuss this conditional law 
as a special case of \cite[Proposition 10 and Theorem 11]{MR2308593}.

\subsection{Organization}

The organization of this paper is as follows. 
In Section \ref{sec: ex} we discuss an example. 
In Section \ref{sec: F-cliques} we prove Proposition \ref{prop: contr W} 
and discuss characterization of stationary probabilities. 
Section \ref{sec: proof} is devoted to the proof of our main theorem, Theorem \ref{thm: main}. 
In Section \ref{sec: non-stationary} we discuss the non-stationary case. 
In Section \ref{sec: app1}, as an appendix, 
we give the proofs for basic facts about the semigroup consisting of mappings. 
Finally, in Section \ref{sec: app2}, we discuss another example 
where the infinite convolution has at least two cluster points.

\subsection*{Acknowledgements}

The authors would like to thank the referee for a lot of valuable comments 
which helped improve the earlier versions of this paper. 
In particular the two appendices are mainly due to the referee.

\section{Example} \label{sec: ex}

Let us investigate an example 
which was discussed in \cite[Subsection 3.3]{MR3023844} 
for mono-particle $ \mu $-evolution. 
We look at it from the viewpoint of multiparticle $ \mu $-evolution. 
See \cite{ISYkkr} for other examples. 

Let $ V = \{ 1,2,3,4,5 \} $. 
We write $ f = [y^1,y^2,y^3,y^4,y^5] $ if $ f:V \to V $ is 
such that $ f1 = y^1,\ldots,f5=y^5 $. 
Consider the two mappings 
\begin{align}
f = [2,3,4,1,5] 
\quad \text{and} \quad 
g = [2,5,5,2,4] . 
\label{}
\end{align}
Let $ \mu = (\delta_f + \delta_g)/2 $ be the uniform law on $ \{ f,g \} $, 
where $ \delta_f $ stands for the Dirac mass at $ f $. 
The marginal transition probability $ P_\mu $ of \eqref{eq: mtp} is given as 
\begin{align}
\pmat{ 
P_\mu(1,\{ 1 \}) & P_\mu(1,\{ 2 \}) & \cdots & P_\mu(1,\{ 5 \}) \\
P_\mu(2,\{ 1 \}) & P_\mu(2,\{ 2 \}) & \cdots & P_\mu(2,\{ 5 \}) \\
\vdots & \vdots & \ddots & \vdots \\
P_\mu(5,\{ 1 \}) & P_\mu(5,\{ 2 \}) & \cdots & P_\mu(5,\{ 5 \}) 
} 
= \frac{1}{2} \pmat{ 
0 & 2 & 0 & 0 & 0 \\
0 & 0 & 1 & 0 & 1 \\
0 & 0 & 0 & 1 & 1 \\
1 & 1 & 0 & 0 & 0 \\
0 & 0 & 0 & 1 & 1 
}. 
\label{eq: Pmu}
\end{align}
Since the fifth power $ P_\mu^5 $ has all positive entries, 
we see that $ P_\mu $ is irreducible aperiodic. 
It is obvious that $ \mu * \lambda $ = $ \lambda $ 
if and only if $ \lambda P_\mu = \lambda $, 
and it is easy to see that 
there exists a unique $ \mu $-invariant probability measure given as 
\begin{align}
\lambda = 
\frac{1}{9} \delta_1 + \frac{2}{9} \delta_2 
+ \frac{1}{9} \delta_3 + \frac{2}{9} \delta_4 
+ \frac{3}{9} \delta_5 . 
\label{}
\end{align}
In \cite[Theorem 1]{MR3023844}, 
for a stationary mono-particle $ \mu $-evolution $ (X,N) $ 
with $ X $ having $ \lambda $ as its common law, 
it was proved that there exists a third noise $ (U_k)_{k \in \bZ} $ such that 
$ \sigma(U_k) \subset \cF^{X,N}_k $ a.s. for $ k \in \bZ $ and 
\begin{align}
\cF^X_k \subset \cF^N_k \vee \sigma(U_k) 
\quad \text{a.s. for $ k \in \bZ $} 
\label{}
\end{align}
with $ \cF^X_{-\infty } $ being trivial a.s. 
and $ \sigma(U_k) $ being independent of $ \cF^N_k $.

Set $ S = \bigcup_{n=1}^{\infty } \{ f,g \}^n $ 
and we would like to apply 
Propositions \ref{prop: IC}, \ref{prop: mapping semigroup} and \ref{prop: contr W}. 
Set 
\begin{align}
e := g^3 = [4,2,2,4,5] . 
\label{}
\end{align}
Note that $ e $ is an idempotent and that 
\begin{align}
efe = e , \quad eg = ge = g , \quad 
g^2 = [5,4,4,5,2] , \quad 
e f^2 = f^2 e = [2,4,4,2,5]. 
\label{}
\end{align}
Since $ f^4 = {\rm id}_V $, the identity of $ V $, 
we see that $ f $ is bijective. 
Note also that 
\begin{align}
g f^2 = f^2 g^2 = [4,5,5,4,2] 
, \quad 
g^2 f^2 = f^2 g = [5,2,2,5,4] . 
\label{}
\end{align}
If we write $ A = \{ 1,3,5 \} $ and $ B = \{ 2,4,5 \} $, 
then $ fB=A $, $ fA = B $ 
and $ gV = gA = gB = A $, 
and hence we see that 
the minimum rank over all mappings of $ S $ is given as $ m_\mu = 3 $. 
We now see that 
the kernel of $ S $ is given as $ K = \{ h \in S : \#(hV) = 3 \} $. 
Since $ K $ is an ideal containing $ g $, we see that 
\begin{align}
S = \{ {\rm id}_V,f,f^2,f^3 \} \cup K . 
\label{}
\end{align}
Since $ e \in E(K) $, we see that $ e $ is a primitive idempotent. 
Let $ K = LGR $ denote the Rees decomposition of $ K $ at $ e $. 

Let us prove that 
\begin{align}
L = \{ e,fe \} 
, \quad 
G = \{ e,g,g^2,ef^2,gf^2,g^2f^2 \} 
, \quad 
R = \{ e,ef \} , 
\label{eq: example LGR}
\end{align}
where 
\begin{align}
fe=[1,3,3,1,5] , \quad 
ef=[2,2,4,4,5] . 
\label{}
\end{align}
We have already seen that $ hV = A $ or $ B $ for all $ h \in K $. 
Let $ \pi_1 = \pi(e) = \{ \{ 1,4 \},\{ 2,3 \},\{ 5 \} \} $ 
and $ \pi_2 = \pi(ef) = \{ \{ 1,2 \},\{ 3,4 \},\{ 5 \} \} $. 
Note that, for every $ s \in S $ and $ k=0,1,2,3 $, 
the partition $ \pi(sgf^k) $ is finer than the partition $ \pi(gf^k) $ 
and hence equal 
since $ \#(sgf^kV) = 3 = \#(gf^kV) $. 
We then see that 
$ \pi(h) = \pi_1 $ or $ \pi_2 $ for all $ h \in K $; 
in fact, $ K $ is contained in $ Sg \cup Sgf \cup Sgf^2 \cup Sgf^3 $, 
$ \pi(g) = \pi(gf^2) = \pi_1 $ and $ \pi(gf) = \pi(gf^3) = \pi_2 $. 

Since $ G $ is the set of all mappings $ h $ of $ S $ such that 
$ \pi(h) = \pi(e) = A $ and $ hV = eV = \pi_1 $, 
it contains the six elements $ e,g,g^2,ef^2,gf^2,g^2f^2 $. 
On the other hand, 
each element of $ G $ induces a permutation on $ A $ and is determined by it, 
so that $ G $ has at most six elements, 
which shows $ G = \{ e,g,g^2,ef^2,gf^2,g^2f^2 \} $. 
Note that $ E(K) $ contains the four elements $ e,fe,ef,fef^3 $, where 
$ fef^3 = (fe)(ef^2)(ef) = [1,1,3,3,5] $. 
On the other hand, 
each element $ h $ of $ E(K) $ is determined by $ hV = A $ or $ B $ 
and $ \pi(h) = \pi_1 $ or $ \pi_2 $, 
so that $ E(K) $ has at most four elements, 
which shows $ E(K) = \{ e,fe,ef,fef^3 \} $. 
We now see that $ L = \{ h \in E(K) : \pi(h)=\pi(e) \} = \{ e,fe \} $ 
and $ R = \{ h \in E(K) : hV = eV \} = \{ e,ef \} $. 
Therefore we obtain \eqref{eq: example LGR}.

Let $ H = e\cS(\eta)e $ be the subgroup of $ G $ in Proposition \ref{prop: IC}. 
Then we have 
$ \mu * \eta^L * \omega_H * \eta^R = \eta^L * \delta_{\gamma} * \omega_H * \eta^R $ 
so that 
$ \mu * \eta^L * \omega_H = \eta^L * \delta_{\gamma} * \omega_H $, 
since $ \eta^R * \delta_e = \delta_e $. 
Let $ \eta^L = \alpha \delta_e + \beta \delta_{fe} $ 
for some $ \alpha,\beta > 0 $ with $ \alpha+\beta = 1 $. 
Since $ gfe = gefe = ge = g $ and $ f^2e = ef^2 $, we have 
\begin{align}
\mu * \eta^L 
= \rbra{ \frac{1}{2} \delta_f + \frac{1}{2} \delta_g } * \rbra{ \alpha \delta_e + \beta \delta_{fe} } 
= \frac{\alpha }{2} \delta_{fe} + \frac{1}{2} \delta_g + \frac{\beta}{2} \delta_{ef^2} . 
\label{eq: compute eta^L}
\end{align}
Since $ fe = (fe) e \in LH $ 
and $ fe = fee \in \cS(\mu * \eta^L * \omega_H) = \cS(\eta^L * \delta_{\gamma} * \omega_H) 
= L \gamma H $, 
we have $ H \cap \gamma H \neq \emptyset $, 
which shows $ H = G $ and we may take $ \gamma = e $. 
We now have 
\begin{align}
\rbra{ \alpha \delta_e + \beta \delta_{fe} } * \omega_G 
=& \eta^L * \omega_G 
= \eta^L * \delta_\gamma * \omega_G 
\label{} \\
=& \mu * \eta^L * \omega_G 
= \rbra{ \frac{\alpha }{2} \delta_{fe} + \frac{1+\beta}{2} \delta_e } * \omega_G , 
\label{}
\end{align}
since $ \delta_e * \omega_G = \delta_g * \omega_G 
= \delta_{ef^2} * \omega_G = \omega_G $. 
Thus $ \alpha = 2/3 $ and $ \beta = 1/3 $, that is, 
\begin{align}
\eta^L = \frac{2}{3} \delta_e + \frac{1}{3} \delta_{fe} . 
\label{}
\end{align}
In the same way we have $ \eta^R = \frac{2}{3} \delta_e + \frac{1}{3} \delta_{ef} $, 
and thus we have obtained that 
\begin{align}
\mu^n \to \eta = \nu = \eta^L * \omega_G * \eta^R . 
\label{}
\end{align}

Note that 
$ fe = [1,3,3,1,5] $ and $ ef = [2,2,4,4,5] $. 
For $ (a,b,c) \in L \times G \times R $, we have 
\begin{align}
\begin{cases}
a=e & \iff \ aV = \{ 2,4,5 \} 
\\
a=fe \!\! & \iff \ aV = \{ 1,3,5 \} 
\end{cases}
, \qquad 
\begin{cases}
c=e & \iff \ c1=c4, \ c2=c3 
\\
c=ef \!\! & \iff \ c1=c2, \ c3=c4 
\end{cases}
. 
\label{}
\end{align}
We note that elements of $ G $ act as permutations over $ \{ 2,4,5 \} $: 
\begin{align}
e(2,4,5) = (2,4,5) 
, \quad 
g(2,4,5) = (5,2,4) 
, \quad 
h(2,4,5) = (4,2,5) . 
\label{}
\end{align}
It is easy to see that 
\begin{align}
W_\mu = 
\{ (x,y,z) : \text{a permutation of $ (2,4,5) $ or $ (1,3,5) $} \} . 
\label{eq: Wmu perm}
\end{align}
We may take a set $ W $ of Proposition \ref{prop: contr W} as 
\begin{align}
W = \{ (2,4,5) \} . 
\label{}
\end{align}
For example, for $ \vx = (3,5,1) \in W_\mu $, 
we see that $ \vx^L = fe $, $ \vx^G = gh $ and $ \vx^W = (2,4,5) $. 

By (iii) of Proposition \ref{prop: contr W}, we see that 
$ \Lambda = \eta^L * \omega_G * \delta_{(2,4,5)} $ 
is the unique $ \mu $-invariant probability measure on $ V^3_\times $. 
Let $ (\bX,N) $ be a stationary tri-particle $ \mu $-evolution 
such that $ \bX $ has $ \Lambda $ as its common law. 
Then we have the factorization 
\begin{align}
\bX_j = \bX^L_j \bX^G_j (2,4,5) 
= \bX_j^L (M^G_{k,j})^{-1} U_k (2,4,5) 
\quad \text{a.s. for $ j \le k $} 
\label{eq: bXj example}
\end{align}
with $ U_k = \bX_k^G $, $ M_j^G = \bX^G_j (\bX^G_{j-1})^{-1} $ 
and $ M^G_{k,j} = \bX^G_k (\bX^G_j)^{-1} =  M^G_k M^G_{k-1} \cdots M^G_{j+1} $, 
and consequently, we obtain the resolution 
\begin{align}
\cF^{\bX}_k = \cG^N_k \vee \sigma(U_k) 
\quad \text{a.s.} 
\quad 
\text{with $ \cG^N_k = \sigma( \bX^L_j,M^G_j : j \le k ) $}
\label{eq: with cGNk}
\end{align}
where the two $ \sigma $-fields 
$ \cF^N_k (\supset \cG^N_k) $ and $ \sigma(U_k) $ are independent. 

We remark that 
\begin{align}
\text{the two processes $ \bX^L $ and $ \bX^G $ are \emph{not} independent.} 
\label{eq: not indep}
\end{align}
In fact, we have 
\begin{align}
& \bP(\bX_0^L = e, \bX_1^L = fe, \bX_0^G = e, \bX_1^G = g ) 
= \bP(\bX_0 = (2,4,5), \bX_1 = feg(2,4,5)) = 0, 
\label{} \\
& \bP(\bX_0^L = e, \bX_1^L = fe) 
= \bP(\bX_0^L=e , N_1 = f) 
= \bP(\bX_0^L=e) \bP(N_1 = f) = \frac{1}{3}, 
\label{} \\
& \bP(\bX_0^G = e, \bX_1^G = g ) 
= \bP(\bX_0^G=e , N_1 = g) 
= \bP(\bX_0^G=e) \bP(N_1 = g) = \frac{1}{12} , 
\label{}
\end{align}
which shows that 
the two events $ \{ \bX_0^L = e, \bX_1^L = fe \} $ and $ \{ \bX_0^G = e, \bX_1^G = g \} $ 
are not independent. 

Note that the first component $ (X^1,N) $ 
is a mono-particle $ \mu $-evolution such that $ X^1 $ has a common law 
\begin{align}
\eta^L * \omega_G * \delta_2 
= \rbra{ \frac{2}{3} \delta_e + \frac{1}{3} \delta_{fe} } * \omega_{\{ 2,4,5 \}} 
= \frac{2}{3} \omega_{\{ 2,4,5 \}} + \frac{1}{3} \omega_{\{ 3,1,5 \}} 
= \lambda , 
\label{eq: etaL omegaG delta2}
\end{align}
where $ \omega_A $ stands for the uniform law on a finite set $ A $. 
Since $ X^1_k = \bX^L_k U_k 2 $ and $ e X^1_k = U_k 2 $, we have 
\begin{align}
\cF^{X^1,N}_k = \cF^N_k \vee \sigma(U_k2) 
\quad \text{a.s. for $ k \in \bZ $}, 
\label{}
\end{align}
where $ U_k2 $ is independent of $ \cF^N_k $. 
We thus conclude that $ (U_k2)_{k \in \bZ} $ is a third noise for $ (X^1,N) $.

\section{F-cliques and stationary probabilities} \label{sec: F-cliques}

Throughout this section 
we suppose all the assumptions of Proposition \ref{prop: contr W} be satisfied.

We borrow several notation from the \emph{road coloring theory}. 
A pair $ \{ x,y \} $ from $ V $ will be called a \emph{deadlock} 
if $ gx \neq gy $ for all $ g \in S := \bigcup_{n=1}^{\infty } \cS(\mu)^n $, 
or in other words, 
$ f_n f_{n-1} \cdots f_1 x \neq f_n f_{n-1} \cdots f_1 y $ 
for all $ n \in \bN $ and $ f_1,\ldots,f_n \in \cS(\mu) $. 
A subset $ F $ of $ V $ will be called an \emph{F-clique} 
if every pair from $ F $ is a deadlock 
and $ F = gV $ for some $ g \in S $. 
Note that a subset $ F $ of $ V $ is an F-clique if and only if 
$ F = gV $ for some $ g \in S $ 
and if, for every $ h \in S $, the restriction $ h|_F $ is injective. 

The F-cliques can be characterized as follows (see also \cite[Lemma 1]{MR3023844}). 

\begin{Lem} \label{lem: F-clique}
For $ g \in S $, 
the set $ gV $ is an F-clique 
if and only if $ \# (gV) = m_{\mu} $. 
Consequently, the image of an F-clique by any mapping of $ S $ is still an F-clique. 
In addition, it holds that 
\begin{align}
\cS(\nu) 
= \{ g \in S : \text{$ gV $ is an F-clique} \} 
= \{ g \in S : \# (gV) = m_{\mu} \} . 
\label{eq: Snu Fcliques}
\end{align}
\end{Lem}

\Proof{
If $ \# (gV) = m_{\mu} $, for any $ f \in S $ we have 
$ m_{\mu} \le \# (fgV) \le \# (gV) = m_{\mu} $ 
so that $ \# (fgV) = m_{\mu} $, 
which implies that $ gV $ is an F-clique. 
Conversely, if $ gV $ is an F-clique, then 
$ \#(fV) \ge \# (fgV) = \# (gV) \ge m_{\mu} $ for any $ f \in S $ 
so that $ \# (gV) = m_{\mu} $. 

Recall that $ \cS(\nu) $ is the kernel, i.e., the unique minimal ideal of 
$ S $ (see (ii) of Proposition \ref{prop: IC}). 
Therefore, to prove \eqref{eq: Snu Fcliques}, 
it suffices to show that 
$ K := \{ g \in S : \# (gV) = m_{\mu} \} $ is a minimal ideal of $ S $. 
It is obvious by definition that $ K $ is an ideal. 
Suppose $ \emptyset \neq IS \cup SI \subset I \subset K $. 
Let $ f \in I $ and $ g \in K $. 
Since $ gf|_{gV}:gV \to gV $ is bijective, 
the mapping $ (gf|_{gV})^r $ is identity for some $ r \in \bN $ 
so that $ (gf)^rg = g $. 
Hence $ g = (gf)^{r-1} gfg \in SIS \subset I $, 
which shows $ I = K $. 
}

By Lemma \ref{lem: F-clique}, 
the set $ W_\mu $ defined in \eqref{eq: Wmu} can be represented as 
\begin{align}
W_\mu = \{ \vx = (x^1,\ldots,x^{m_\mu}) 
: \text{$ \{ x^1,\ldots,x^{m_\mu} \} $ is an F-clique} \} . 
\label{eq: Wmu rep}
\end{align}

\begin{Lem} \label{lem1} 
For any $ \vx,\vx' \in eW_\mu $, 
the two sets $ G \vx $ and $ G \vx' $ 
are either equal or disjoint. 
\end{Lem}

\Proof{
Suppose $ G \vx $ 
and $ G \vx' $ have a common element 
$ g\vx = g'\vx' $ for some $ g,g' \in G $. 
We then obtain that 
\begin{align}
G \vx = G g \vx 
= G g' \vx' 
= G \vx'. 
\label{}
\end{align}
The proof is complete. 
}

We now prove Proposition \ref{prop: contr W}. 

\Proof[Proof of Proposition \ref{prop: contr W}]{
(i) 
By \eqref{eq: Snu Fcliques} and \eqref{eq: Wmu rep}, 
we have $ W_\mu = \cS(\nu) W_\mu = LGRW_\mu $. Let us prove that 
\begin{align}
RW_\mu = GW_\mu = eW_\mu 
\quad \text{and} \quad 
LW_\mu = W_\mu . 
\label{eq: RGLWmu}
\end{align}
Since $ eW_\mu \subset RW_\mu \subset eLGRW_\mu = eW_\mu $ 
and $ eW_\mu \subset GW_\mu \subset eLGRW_\mu = eW_\mu $, 
we have $ RW_\mu = GW_\mu = eW_\mu $. 
Hence $ LW_\mu = LLGRW_\mu = LGRW_\mu = W_\mu $. 
We now obtain \eqref{eq: RGLWmu}.

We thus obtain 
\begin{align}
W_\mu = LGRW_\mu = LGGW = LGW . 
\label{eq: Wmu LGW}
\end{align}

(ii) 
Note that 
\begin{align}
e W_\mu = GW = \bigcup_{\vw \in W} G \vw . 
\label{}
\end{align}
By the minimality of $ W $, we see that the sets $ G\vw $ for $ \vw \in W $ are disjoint; 
in fact, if $ G\vw $ and $ G\vw' $ are not disjoint 
with some distinct elements $ \vw $ and $ \vw' $ of $ W $, 
then $ G\vw = G\vw' $ by Lemma \ref{lem1} and so $ eW_\mu = G(W \setminus \{ \vw' \}) $, 
which contradicts the minimality of $ W $. 

We have only to prove injectivity of the product 
$ L \times G \times W \ni (f,g,\vw) \mapsto fg\vw \in W_\mu $. 
Suppose $ fg\vw = f'g'\vw' $. 
Since $ eL = \{ e \} $, we have $ g\vw = g'\vw' $, which implies $ G\vw = G\vw' $. 
By the above argument, we have $ \vw = \vw' $. 

If we write $ \vw = (w^1,\ldots,w^{m_{\mu}}) $, 
then $ \{ w^1,\ldots,w^{m_{\mu}} \} = eV $, because 
$ ( w^1,\ldots,w^{m_\mu} ) \in W \subset eW_\mu $ so that 
the points $ w^1,\ldots,w^{m_\mu} $ are distinct elements of $ eV $, 
and $ \#(eV) = m_\mu $ by \eqref{eq: Snu Fcliques}. 
Hence the identity $ fg\vw = f'g'\vw $ implies that $ fg=f'g' $ on $ eV $. 
Since $ g=ge $ and $ g'=g'e $, 
we see that $ fg=f'g' $ on $ V $, which implies $ f=f' $ and $ g=g' $. 

(iii) 
Let $ \Lambda_W \in \cP(W) $ and set $ \Lambda = \eta^L * \omega_G * \Lambda_W $. 
Since $ Re = \{ e \} $, we have 
$ \nu * \delta_e = \eta^L * \omega_G * \eta^R * \delta_e 
= \eta^L * \omega_G * \delta_e = \eta^L * \omega_G $. 
We thus see that 
\begin{align}
\mu * \Lambda 
= \mu * \nu * \delta_e * \Lambda_W 
= \nu * \delta_e * \Lambda_W 
= \eta^L * \omega_G * \Lambda_W 
= \Lambda . 
\label{eq: mu*Lambda = Lambda}
\end{align}

Conversely, suppose $ \Lambda \in \cP(V^{m_\mu}_\times) $ be $ \mu $-invariant. 
Since $ \Lambda = \mu * \Lambda $, 
we have $ \nu * \Lambda = \lim_{n \to \infty } 
\frac{1}{n} \sum_{k=1}^n \mu^k * \Lambda = \Lambda $, 
and hence $ \Lambda = \nu * \Lambda = \eta^L * \omega_G * \eta^R * \Lambda $. 
Since $ \cS(\nu) \cS(\Lambda) = \cS(\nu * \Lambda) = \cS(\Lambda) \subset V^{m_\mu}_\times $ 
and by \eqref{eq: Wmu rep}, 
we have $ \cS(\Lambda) \subset W_\mu $. 
Let us explain it in detail. If 
$ (x_1,\ldots,x_{m_\mu}) \in \cS(\Lambda) $, then $ x_1,\ldots,x_{m_\mu} $ are distinct. 
Moreover, since $ \cS(\Lambda) = \cS(\nu) \cS(\Lambda) $, 
we have $ (x_1,\ldots,x_{m_\mu}) = k(y_1,\ldots,y_{m_\mu}) $ for some $ k \in \cS(\nu) $ and
$ (y_1,\ldots,y_{m_\mu}) \in \cS(\Lambda) \subset V^{m_\mu}_\times $. 
Thus $ \{ x_1,\ldots,x_{m_\mu} \} \subset kV $, so $ \{ x_1,\ldots,x_{m_\mu} \} = kV $ 
since $ \#(kV)=m_\mu $. By Lemma \ref{lem: F-clique}, 
$ \{ x_1,\ldots,x_{m_\mu} \} $ is an F-clique, i.e.
$ (x_1,\ldots,x_{m_\mu}) \in W_\mu $ by \eqref{eq: Wmu rep}. 

Since $ \cS(\eta^R * \eta^L) \subset RL \subset G $ and thus 
$ \eta^R * \Lambda = (\eta^R * \eta^L) * \omega_G * \eta^R * \Lambda 
= \omega_G * \eta^R * \Lambda $, we have 
$ \cS(\eta^R * \Lambda) = \cS(\omega_G * \eta^R * \Lambda) 
\subset GRW_\mu = eW_\mu = GW $. Hence 
\begin{align}
\Lambda = (\eta^L * \omega_G) * (\eta^R * \Lambda) 
=& \sum_{\vx \in GW} (\eta^R * \Lambda)\{ \vx \} \, (\eta^L * \omega_G * \delta_{\vx}) 
\label{} \\
=& \sum_{\vx \in GW} (\eta^R * \Lambda)\{ \vx \} \, (\eta^L * \omega_G * \delta_{\vx^W}) 
= \eta^L * \omega_G * \Lambda_W , 
\label{}
\end{align}
where we take 
\begin{align}
\Lambda_W = \sum_{\vx \in GW} (\eta^R * \Lambda) \{ \vx \} \delta_{\vx^W} . 
\label{}
\end{align}
The proof is now complete. 
}

Let us prove Corollary \ref{cor: extension}. 

\Proof[Proof of Corollary \ref{cor: extension}]{
Let $ \lambda \in \cP(V) $ be $ \mu $-invariant. 
We then have 
\begin{align}
\lambda 
= \lim_{n \to \infty } \frac{1}{n} \sum_{k=1}^n \mu^k * \lambda 
= \nu * \lambda = \eta^L * \omega_G * \eta^R * \lambda . 
\label{eq: lambda Cesaro}
\end{align}
For $ w \in V $, we write 
\begin{align}
n(w) := \# (\{ (x^1,\ldots,x^{m_\mu}) \in W : x^1 = w \}) , 
\label{}
\end{align}
which turns out to be positive by definition of $ W $. 
We then define 
\begin{align}
\Lambda_W \{ (v^1,\ldots,v^{m_\mu}) \} 
= \frac{1}{n(v^1)} (\eta^R * \lambda) \{ v^1 \} 
, \quad (v^1,\ldots,v^{m_\mu}) \in W 
\label{eq: LambdaW defi}
\end{align}
so that $ \Lambda_W \in \cP(W) $, 
and define $ \Lambda = \eta^L * \omega_G * \Lambda_W $. 
By (iii) of Proposition \ref{prop: contr W}, we see that $ \Lambda $ is $ \mu $-invariant. 
Let us compute the marginal of $ \Lambda $ in the first coordinate. 
For any $ v \in V $, we have 
\begin{align}
& \Lambda \{ (x^1,\ldots,x^{m_\mu}) \in V^{m_\mu}_\times : x^1 = v \} 
\label{} \\
=& \sum_{\begin{subarray}{c} f \in L \\ g \in G \end{subarray}} \, 
\sum_{(v^1,\ldots,v^{m_\mu}) \in W} 
\eta^L \{ f \} \omega_G \{ g \} 
\Lambda_W \{ (v^1,\ldots,v^{m_\mu}) \} 1_{\{ fgv^1 = v \}} 
\label{} \\
=& \sum_{\begin{subarray}{c} f \in L \\ g \in G \end{subarray}} \, 
\sum_{(v^1,\ldots,v^{m_\mu}) \in W} 
\eta^L \{ f \} \omega_G \{ g \} 
\cdot \frac{1}{n(v^1)} (\eta^R * \lambda) \{ v^1 \} 
1_{\{ fgv^1 = v \}} 
\label{} \\
=& \sum_{\begin{subarray}{c} f \in L \\ g \in G \end{subarray}} \, \sum_{v^0 \in V} 
\eta^L \{ f \} \omega_G \{ g \} 
(\eta^R * \lambda) \{ v^0 \} 
1_{\{ fgv^0 = v \}} 
\cdot \frac{1}{n(v^0)} \sum_{(v^1,\ldots,v^{m_\mu}) \in W} 
1_{\{ v^1 = v^0 \}} 
\label{eq: v0 sum} \\
=& (\eta^L * \omega_G * \eta^R * \lambda) \{ v \} = \lambda \{ v \} . 
\label{}
\end{align}
The proof is complete. 
}

\section{Proof of our main theorem} \label{sec: proof}

Throughout this section 
we suppose all the assumptions of Theorem \ref{thm: main} be satisfied. 
We divide the proof of Theorem \ref{thm: main} into several steps.

\subsection{Factorizing $ \bX_k $ into $ LG $- and $ W $-factors} \label{sec: Fac bXk LG W}

For any fixed $ k \in \bZ $, we see, by Proposition \ref{prop: contr W}, that 
$ \bX_k \in LGW $ a.s. and the law of $ \bX_k $ 
is equal to $ \eta^L * \omega_G * \Lambda_W $ for some $ \Lambda_W \in \cP(W) $, 
which shows that the three random variables $ \bX_k^L $, $ \bX_k^G $ and $ \bX_k^W $ 
are independent and for the marginal laws we have $ \bX_k^L \dist \eta^L $, 
$ \bX_k^G \dist \omega_G $ and $ \bX_k^W \dist \Lambda_W $. 
Hence we have shown Claim (i) of Theorem \ref{thm: main}. 

Let us focus on the factor $ \bX_k^L \bX_k^G $ 
in the factorization $ \bX_k = \bX_k^L \bX_k^G \bX_k^W $ for $ k \in \bZ $. 

\begin{Prop} \label{prop: Yk}
Set $ Y_k = \bX_k^L \bX_k^G $ for $ k \in \bZ $ and $ Y = (Y_k)_{k \in \bZ} $. 
Then the following hold: 
\begin{enumerate}

\item 
$ (Y,N) $ is a $ \mu $-evolution 
such that the sequence $ Y $ has a common law $ \eta^L * \omega_G $. 

\item 
There exists a $ W $-valued random variable $ \bZ_W $ such that 
$ \bX_k^W = \bZ_W $ a.s. for $ k \in \bZ $. 

\item 
$ (Y,N) $ and $ \bZ_W $ are independent. 

\end{enumerate}
\end{Prop}

\Proof{
By the argument in the beginning of this subsection, 
we see that, for every $ k \in \bZ $, the law of $ Y_k = \bX_k^L \bX_k^G $ 
is $ \eta^L * \omega_G $. 

Note that 
\begin{align}
Y_k \bX_k^W = \bX_k = N_k \bX_{k-1} = (N_k Y_{k-1}) \bX_{k-1}^W 
\quad \text{a.s.} 
\label{eq: Yk bXkW}
\end{align}
Since $ Y_{k-1} \in LG $ and 
$ N_k Y_{k-1} \in SLG = SSe \subset Se = LG $ by Proposition \ref{prop: kernel}, 
we see, by Proposition \ref{prop: contr W}, that 
\begin{align}
Y_k = N_k Y_{k-1} 
\quad \text{and} \quad 
\bX_k^W = \bX_{k-1}^W 
\quad \text{a.s.} 
\label{}
\end{align}
Since $ N_k $ is independent of $ \cF^Y_{k-1} (\subset \cF^{\bX}_{k-1}) $, 
we see that $ (Y,N) $ is a $ \mu $-evolution. 
We now obtain Claims (i) and (ii) 
(and consequently we have shown Claim (iii) of Theorem \ref{thm: main}). 

Let $ k \in \bZ $ be fixed. 
By the above argument, we see that $ Y_k = \bX_k^L \bX_k^G $ is independent of $ \bZ_W $. 
Since $ \{ N_j : j>k \} $ is independent of $ \{ Y_k,\bZ_W \} $ 
and since $ Y_j = N_j N_{j-1} \cdots N_{k+1} Y_k $ for $ j>k $, 
we see that $ \{ (Y_j,N_j) : j>k \} $ is independent of $ \bZ_W $. 
Since $ k \in \bZ $ is arbitrary, we obtain Claim (iii). 
The proof is complete. 
}

\subsection{Factorizing $ \bX^G_k $ into $ C $- and $ H $-factors} \label{sec: Fac bXG CH}

For $ f \in \cS(\nu) = LGR $, we write 
$ f^C = (f^G)^C $ and $ f^H = (f^G)^H $. 
Consequently, the mapping $ z \mapsto (z^L,z^C,z^H,z^R) $ 
is the inverse of the product mapping 
$ L \times C \times H \times R \ni (x,c,h,y) \mapsto x c h y \in \cS(\nu) $. 
For $ \vx \in LGW $, we write $ \vx^C = (\vx^G)^C $ and $ \vx^H = (\vx^G)^H $. 
Consequently, the mapping $ \vx \mapsto (\vx^L,\vx^C,\vx^H,\vx^W) $ 
is the inverse of the product mapping 
$ L \times C \times H \times W \ni (x,c,h,\vw) \mapsto x c h \vw \in LGW $. 

Since $ H $ is a normal subgroup of $ G $, we have 
\begin{align}
(g_1 g_2)^C H 
= (g_1 g_2) H 
= (g_1 H) (g_2 H) 
= (g_1^C H) (g_2^C H) 
= (g_1^C g_2^C) H , 
\label{}
\end{align}
so that $ (g_1 g_2)^C = (g_1^C g_2^C)^C $.

We proceed to prove part of Theorem \ref{thm: main}. 

\begin{Prop} \label{prop: NGk}
Claim \eqref{eq: cGNk cFNk} of Theorem \ref{thm: main} holds 
and it holds that 
\begin{align}
\text{$ \bX_k^C H = \gamma^k Y_C H $ a.s. for $ k \in \bZ $ 
for some $ C $-valued random variable $ Y_C $.}
\label{eq: (iii) former half}
\end{align}
Consequently, 
$ \bX_k^C = (\gamma^k Y_C)^C $ a.s. for $ k \in \bZ $. 
\end{Prop}

\Proof{
Set 
\begin{align}
N_{k,l} := N_k N_{k-1} \cdots N_{l+1} 
, \quad k > l . 
\label{}
\end{align}
Since $ e \in S = \bigcup_{n=1}^{\infty } \cS(\mu)^n $, 
we can find $ f_1,f_2,\ldots,f_{n_0} \in \cS(\mu) $ such that 
$ f_{n_0} f_{n_0-1} \cdots f_1 = e $, and hence we have 
\begin{align}
T^e_k := \sup \{ l<k-n_0 : N_{l+n_0,l} = e \} 
> - \infty \quad \text{a.s.} 
\label{eq: Tek}
\end{align}
We now see that 
$ N_{k,T^e_k} = N_{k,T^e_k+n_0} N_{T^e_k+n_0,T^e_k} = N_{k,T^e_k+n_0} e \in Se = LG $ 
by Proposition \ref{prop: kernel}. 

Let us prove Claim \eqref{eq: cGNk cFNk}. 
Since $ \bX_k = N_{k,T^e_k} \bX_{T^e_k} $, we have 
\begin{align}
\bX_k^L = (N_{k,T^e_k})^L 
, \quad 
\bX_k^G = (N_{k,T^e_k})^G \, \bX_{T^e_k}^G 
\quad \text{a.s.} 
\label{}
\end{align}
Hence we obtain $ \bX_k^L \in \cF^N_k $ a.s. 
Since $ \bX_k = N_k \bX_{k-1}^L \bX_{k-1}^G \bX_{k-1}^W $ 
and $ N_k \bX^L_{k-1} \in SL = SLe \subset Se = LG $ 
by Proposition \ref{prop: kernel}, we have 
\begin{align}
\bX_k^L = (N_k \bX_{k-1}^L)^L 
, \quad 
\bX_k^G = (N_k \bX_{k-1}^L)^G \, \bX_{k-1}^G 
\quad \text{a.s.} 
\label{eq: 4.8}
\end{align}
Hence we obtain $ M_k^G = \bX_k^G (\bX_{k-1}^G)^{-1} = (N_k \bX_{k-1}^L)^G \in \cF^N_k $ a.s. 
We thus obtain Claim \eqref{eq: cGNk cFNk}. 

Let $ \xi $ be a random variable such that 
$ \xi \dist \omega_H $ and $ \xi $ is independent of $ (\bX,N) $. 
Let $ k \in \bZ $. 
By $ N_k \bX_{k-1}^L \in LG $, we have 
\begin{align}
M_k^G \xi = (N_k \bX_{k-1}^L)^G \xi 
= (N_k \bX_{k-1}^L \xi)^G . 
\label{}
\end{align}
Since $ \omega_H * \eta^R * \delta_e = \omega_H $ by $ Re = \{ e \} $ and $ e \in H $, 
we have $ \eta * \delta_e = \eta_L * \omega_H $, so 
\begin{align}
N_k \bX_{k-1}^L \xi 
\dist \mu * \eta^L * \omega_H 
= \mu * \eta * \delta_e = \eta^L * \delta_{\gamma} * \omega_H * \eta^R * \delta_e 
= \eta^L * \delta_{\gamma} * \omega_H , 
\label{eq: NkbXk-1Lxi}
\end{align}
we have $ M_k^G \xi \dist \delta_{\gamma} * \omega_H \dist \gamma \xi $, 
which shows $ (M_k^G)^C = \gamma $ a.s. for $ k \in \bZ $. 
We now see that 
\begin{align}
\bX_k^C H = \bX_k^G H 
= M_k^G \bX_{k-1}^G H 
= (M_k^G)^C (\bX_{k-1}^G)^C H 
= \gamma \bX_{k-1}^C H 
\quad \text{a.s. for $ k \in \bZ $}. 
\label{}
\end{align}
Hence $ \gamma^{-k} \bX_k^C H = \gamma^{-(k-1)} \bX_{k-1}^C H $ a.s., 
which shows that there exists a $ C $-valued random variable $ Y_C $ such that 
$ \gamma^{-k} \bX_k^C H = Y_C H $ a.s. for all $ k \in \bZ $, 
which yields \eqref{eq: (iii) former half}. 
The proof is now complete. 
}

\begin{Rem} \label{rem: BL}
Let $ (\bX,N) $ be a stationary $ m_\mu $-particle $ \mu $-evolution. 
Let us pick up \cite[Proposition 10 and Theorem 11]{MR2308593}, 
whose conclusion applied to our $ \mu $-evolutions is as follows: 
\emph{If $ \bX $ is an irreducible aperiodic recurrent Markov chain on $ W_{\mu} $, 
then the random set 
\begin{align}
R_0 := \bigcap_{j<0} N_{0,j} W_{\mu} 
\label{eq: R0bigcap}
\end{align}
has exactly $ m_\mu $ elements, 
the remote past noise $ \cF^{\bX}_{-\infty } $ is trivial, 
and the conditional law of $ \bX_0 $ given $ \sigma(N) $ is uniform on $ R_0 $.} 
Let us derive this result from our results. 
Since $ f \in eS $ implies $ fV = eV $ by (v) of Proposition \ref{prop: mapping semigroup}, 
and since $ eW_\mu $ is a subset of the set of all permutations of $ eV $ by \eqref{eq: Wmu rep}, 
we see that $ f \in eS $ implies $ fW_\mu = eW_\mu $. 
For $ j < T_0^e $, we have 
\begin{align}
N_{0,j} W_\mu 
= N_{0,T_0^e} N_{T_0^e,j} W_\mu 
= (N_{0,T_0^e})^L (N_{0,T_0^e})^G N_{T_0^e,j} W_\mu 
= \bX^L_0 eW_\mu , 
\label{}
\end{align}
since $ (N_{0,T_0^e})^G N_{T_0^e,j} $ takes values in $ GS = eGS \subset eS $. 
Hence we obtain 
\begin{align}
R_0 = \bX^L_0 eW_\mu . 
\label{eq: R0bXL0eWmu}
\end{align}
(We did not need aperiodicity so far.) 
Let us now consider the special case where 
$ \bX $ is an irreducible aperiodic Markov chain on $ W_{\mu} $. 
For every $ \vw \in W $, 
the law $ \mu^n * \delta_{\vw} $ converges 
to the unique $ \mu $-invariant probability measure $ \Lambda \in \cP(W_\mu) $. 
This shows that $ \Lambda = \eta * \delta_{\vw} = \eta^L * \omega_H * \delta_{\vw} $. 
By (iii) of Proposition \ref{prop: contr W}, we obtain 
$ H = G $ and $ W = \{ \vw \} $ is a singleton, 
and consequently $ \cF^{\bX}_{-\infty } $ is trivial, 
and $ eW_\mu = GW = G \{ \vw \} $. 
As we have seen it in Section \ref{sec: HR}, the conditional law 
$ \bP(\bX_0 \in \cdot \mid \sigma(N)) $ is uniform on the random finite set 
$ \{ \bX^L_0 g \vw : g \in G \} = \bX^L_0 eW_\mu = R_0 $. 
\end{Rem}

\subsection{Finding the third noise}

The following lemma plays a key role. 

\begin{Lem} \label{lem: key} 
Let $ \{ a_n \} $ and $ \{ b_n \} $ be 
two deterministic sequences of $ \cS(\nu) $. Then 
\begin{align}
(\delta_{a_n} * \mu^n * \delta_{b_n})^H \cdist \omega_H . 
\label{eq: fn N1n hn}
\end{align}
\end{Lem}

\Proof{
Since $ \cP(H) $ is compact, it suffices to show that 
$ \omega_H $ is the only one cluster point 
of the sequence $ \{ (\delta_{a_n} * \mu^n * \delta_{b_n})^H \} $. 
Let $ \{ n(m) \} $ be a subsequence of $ \bN $ such that 
$ (\delta_{a_n} * \mu^n * \delta_{b_n})^H $ converges as $ n \to \infty $. 
Taking a further subsequence if necessary, 
we can and do assume that 
$ a_{n(m)} \to a_0 $ and $ b_{n(m)} \to b_0 $ for some $ a_0,b_0 \in \cS(\nu) $ and 
\begin{align}
\mu^{n(m)} \to \mu^r * \eta = \eta^L * \delta_{\gamma^r} * \omega_H * \eta^R 
\label{}
\end{align}
for some $ r = 0,1,\ldots,p-1 $. Hence we have 
\begin{align}
(\delta_{a_{n(m)}} * \mu^{n(m)} * \delta_{b_{n(m)}})^H 
\cdist (\delta_{a_0} * \eta^L * \delta_{\gamma^r} * \omega_H * \eta^R * \delta_{b_0})^H . 
\label{}
\end{align}
Since $ RL \subset H $ and $ g^{-1} H g = H $ for all $ g \in G $, we have 
\begin{align}
& (\delta_{a_0} * \eta^L * \delta_{\gamma^r} * \omega_H * \eta^R * \delta_{b_0})^G 
\label{} \\
=& \delta_{a^G_0} * \delta_{a^R_0} * \eta^L 
* \delta_{\gamma^r} * \omega_H * \eta^R 
* \delta_{b^L_0} * \delta_{b^G_0} 
\label{} \\
=& \delta_{a_0^C \gamma^r} * \delta_{\gamma^{-r} a^H_0 \gamma^r} 
* (\delta_{\gamma^{-r}} * \delta_{a^R_0} * \eta^L * \delta_{\gamma^r}) * \omega_H * (\eta^R * \delta_{b_0^L}) * \delta_{(b_0^C b_0^H (b_0^C)^{-1})} * \delta_{b_0^C} 
\label{} \\
=& \delta_{a_0^C \gamma^r} * \omega_H * \delta_{b_0^C} 
= \delta_{a_0^C \gamma^r b_0^C} * \omega_H 
= \delta_{(a_0^C \gamma^r b_0^C)^C} * \omega_H , 
\label{}
\end{align}
which yields $ (\delta_{a_0} * \eta^L * \delta_{\gamma^r} * \omega_H * \eta^R * \delta_{b_0})^H 
= \omega_H $. 
We thus obtain \eqref{eq: fn N1n hn}. 
}

We proceed to prove part of Theorem \ref{thm: main}. 

\begin{Prop} \label{prop: third noise main}
Let $ k \in \bZ $ be fixed and set $ U_k := \bX^H_k = Y^H_k $. 
Then $ U_k \dist \omega_H $ 
and the three $ \sigma $-fields $ \sigma(N) $, $ \cF^{\bX}_{-\infty } $ and 
$ \sigma(U_k) $ are independent. 
(Consequently Claims (ii) and \eqref{eq: main indep} of Theorem \ref{thm: main} hold.) 
\end{Prop}

\Proof{
Set 
\begin{align}
\cF^N_{k,l} = \sigma( N_k,N_{k-1},\ldots,N_{l+1} ) 
, \quad k>l . 
\label{}
\end{align}
Let $ k \in \bZ $ be fixed and let $ \varphi:H \to \bR $ be a test function. 
Let $ l < k $, $ k_0 > k $, $ n \in \bN $, $ A \in \cF^N_{k_0,l} $ 
and $ B \in \cF^{\bX}_{-\infty } $. 
Recall that the symbols $ T^e_l $ and $ N_{k,l} $ 
have been introduced in the proof of Proposition \ref{prop: NGk}. 
Since $ T^e_l - n < T^e_l < l $, we see that 
the three $ \sigma $-fields 
$ \sigma( N_{k,T^e_l},1_A ) $, $ \sigma( N_{T^e_l,T^e_l - n} ) $ 
and $ \sigma( Y_{T^e_l-n} , 1_B ) $ are independent. 
Since $ N_{T^e_l,T^e_l - n} \dist \mu^n $, we have 
\begin{align}
\bE \sbra{ \varphi(U_k) 1_A 1_B } 
=& \bE \sbra{ \varphi(Y_k^H) 1_A 1_B } 
\label{} \\
=& \bE \sbra{ \varphi \rbra{ (N_{k,T^e_l} N_{T^e_l,T^e_l - n} Y_{T^e_l-n})^H } 1_A 1_B } 
\label{} \\
=& \bE \sbra{ \left. \int \varphi \, \d (\delta_{a} * \mu^n * \delta_{b_n})^H 
\right|_{\begin{subarray}{l} a = N_{k,T^e_l} \\ b_n = Y_{T^e_l-n} \end{subarray}} 
1_A 1_B } . 
\label{eq: bE'}
\end{align}
Noting that $ N_{k,T^e_l} \in \cS(\nu) $ (see the proof of Proposition \ref{prop: NGk}) 
and $ Y_{T^e_l-n} \in \cS(\nu) $, 
we apply Lemma \ref{lem: key} to see that 
\begin{align}
\eqref{eq: bE'} 
\tend{}{n \to \infty } 
\int \varphi \, \d \omega_H \cdot \bE \sbra{ 1_A 1_B } 
= \int \varphi \, \d \omega_H \cdot \bP(A) \bP(B) . 
\label{}
\end{align}
Since $ l < k $ and $ k_0 > k $ are arbitrary, we obtain 
\begin{align}
\bE \sbra{ \varphi(U_k) 1_A 1_B } 
= \int \varphi \, \d \omega_H \cdot \bP(A) \bP(B) 
, \quad A \in \sigma(N) , \ B \in \cF^{\bX}_{- \infty } , 
\label{eq: varphiU_k 1_A 1_B}
\end{align}
which leads to the desired result. 
}

\subsection{Determining the remote past noise}

We need the following lemma. 

\begin{Lem} \label{lem: indep}
Let $ (\Omega,\cF,\bP) $ be a probability space 
and let $ \cA $, $ \cB $ and $ \cC $ be three sub-$ \sigma $-fields of $ \cF $. 
Suppose that $ \cA \subset \cB \vee \cC $ a.s. 
and that $ \cA \vee \cB $ be independent of $ \cC $. 
Then $ \cA \subset \cB $ a.s. 
\end{Lem}

The proof of Lemma \ref{lem: indep} can be found in \cite[Section 2.2]{CY}, 
and so we omit it. 

We shall now complete the proof of Theorem \ref{thm: main}. 

\Proof[Proof of Theorem \ref{thm: main}]{
What remains unproved are Claims (iv), (v) and (vi). 

We have shown that $ \bX_k^C = (\gamma^k Y_C)^C $, $ \bX_k^H = U_k $ and $ \bX_k^W = \bZ_W $. 
Let $ j \le k $. 
Since $ \bX_k^G = M_k^G \bX_{k-1}^G $ by definition of $ M_k^G $, 
we have $ \bX_k^G = M_{k,j}^G \bX_j^G $. 
Hence we obtain 
\begin{align}
\bX_j 
= \bX_j^L \bX_j^G \bZ_W 
= \bX_j^L (M_{k,j}^G)^{-1} \bX_k^G \bZ_W 
= \bX_j^L (M_{k,j}^G)^{-1} (\gamma^k Y_C)^C U_k \bZ_W 
\quad \text{a.s.}, 
\label{}
\end{align}
which shows Claim (iv) and leads to 
\begin{align}
\cF^{\bX}_k = \cG^N_k \vee \sigma( Y_C,\bZ_W ) \vee \sigma(U_k) 
\quad \text{a.s.} 
\label{eq: resol-}
\end{align}
Since $ \sigma(Y_C,\bZ_W) \subset \cF^{\bX}_{-\infty } $ a.s., 
which is obvious by definition, 
and by \eqref{eq: main indep}, 
we can apply Lemma \ref{lem: indep} 
for $ \cA = \cF^{\bX}_{-\infty } $, $ \cB = \sigma(Y_C,\bZ_W) $ 
and $ \cC = \cF^N_k \vee \sigma(U_k) $, 
and hence we obtain 
$ \cF^{\bX}_{-\infty } \subset \sigma(Y_C,\bZ_W) $ a.s. 
We thus obtain \eqref{eq: bX-infty}. 
Combining \eqref{eq: resol-} and \eqref{eq: bX-infty}, 
we obtain \eqref{eq: resol}, 
which shows Claim (v). 

Since $ \bX_k = \bX_k^L (\gamma^k Y_C)^C U_k \bZ_W $ 
and since $ \Lambda = \eta^L * \omega_G * \Lambda_W $, 
we see that 
$ (\gamma^k Y_C)^C U_k $ and $ \bZ_W $ are independent 
and that $ (\gamma^k Y_C)^C U_k \dist \omega_G $ 
and $ \bZ_W \dist \Lambda_W $. 
Since $ \omega_G = \omega_C * \omega_H $, 
we see that $ \gamma^k Y_C H $ and $ U_k $ are independent, 
$ \gamma^k Y_C H \dist \omega_{G/H} $ and $ U_k \dist \omega_H $, 
which yields that 
$ Y_C $ and $ U_k $ are independent 
and $ Y_C \dist \omega_C $. 
We now obtain Claim (vi). 

The proof of Theorem \ref{thm: main} is therefore complete. 
}

\section{The non-stationary case} \label{sec: non-stationary}

Throughout this section we adopt the settings of Subsection \ref{subsec: main results}. 

\begin{Prop} \label{prop: num of distinct}
Let $ (\bX,N) $ be an $ m $-particle $ \mu $-evolution. 
Then, for almost every sample path 
and for any $ 1 \le i < j \le m $, 
either [$ X^i = X^j $ (i.e., $ X^i_k = X^j_k $ holds for all $ k \in \bZ $)] 
or [the pair $ \{ X^i_k,X^j_k \} $ is a deadlock for all $ k \in \bZ $]. 
Consequently, the number of distinct states among $ \{ X^1_k,\ldots,X^m_k \} $ 
does not depend upon $ k \in \bZ $ a.s. 
\end{Prop}

\Proof{
Recall $ T^e_k $, which has been introduced in the proof of Proposition \ref{prop: NGk}. 
We then have, for almost every sample path, 
$ \bX_k = N_{k,T^e_k} \bX_{T^e_k} \in LGV^m $ for all $ k \in \bZ $. 
By Lemma \ref{lem: F-clique}, we see that 
[every distinct pair from $ \{ X^1_k,\ldots,X^m_k \} $ is a deadlock] for all $ k \in \bZ $. 

Suppose that $ X^i \neq X^j $. 
Then $ X^i_{k_0} \neq X^j_{k_0} $ for some $ k_0 \in \bZ $. 
Since the pair $ \{ X^i_{k_0},X^j_{k_0} \} $ is a deadlock, 
we have $ X^i_k \neq X^j_k $ for all $ k \ge k_0 $. 
We also see that $ X^i_k \neq X^j_k $ for all $ k \le k_0 $, 
since $ (X^i_{k_0},X^j_{k_0}) = N_{k_0,k} (X^i_k,X^j_k) $. 
For every $ k \in \bZ $, we see that $ \{ X^i_k,X^j_k \} $ is a distinct pair, 
and hence is a deadlock. 
The proof is now complete. 
}

\begin{Prop} \label{prop: char Lambda k}
For a sequence $ (\Lambda_k)_{k \in \bZ} $ from $ \cP(V^{m_\mu}_\times) $, 
the following are equivalent: 
\begin{enumerate}

\item 
$ \Lambda_k = \mu * \Lambda_{k-1} $ for all $ k \in \bZ $. 

\item 
There exist $ \Lambda^0_W,\ldots,\Lambda^{p-1}_W \in \cP(W) $ 
and constants $ \alpha_0,\ldots,\alpha_{p-1} \ge 0 $ 
with $ \alpha_0 + \cdots + \alpha_{p-1} = 1 $ 
such that 
\begin{align}
\Lambda_k = \sum_{r=0}^{p-1} \alpha_r \eta^L * \delta_{\gamma^{k+r}} * \omega_H * \Lambda^r_W 
\quad \text{for all $ k \in \bZ $}. 
\label{eq: Lambdak}
\end{align}
(Consequently, $ \Lambda_{k+p} = \Lambda_k $ for all $ k \in \bZ $, since $ \gamma^p \in H $.) 
\end{enumerate}
\end{Prop}

\Proof{
Suppose that Condition (ii) is satisfied. 
Since 
\begin{align}
\mu * (\eta^L * \delta_{\gamma^{k+r}} * \omega_H) 
= \eta^L * \delta_{\gamma^{k+1+r}} * \omega_H , 
\label{eq: mu*...}
\end{align}
we obtain $ \mu * \Lambda_{k-1} = \Lambda_k $ for all $ k \in \bZ $, 
which shows Condition (i). 
Since $ \gamma^p \in H $, we have $ \delta_{\gamma^p} * \omega_H = \omega_H $, 
which yields $ \Lambda_{k+p} = \Lambda_k $ for all $ k \in \bZ $. 

Suppose that Condition (i) is satisfied. 
Iterating the relation $ \Lambda_k = \mu * \Lambda_{k-1} $, 
we have $ \Lambda_k = \mu^p * \Lambda_{k-p} $ for all $ k \in \bZ $. 
Let $ k \in \bZ $ be fixed. 
Since $ \cP(V^{m_\mu}_\times) $ is compact, 
we have $ \Lambda_{-pn(m)} \to \Lambda_* $ 
for some subsequence $ \{ n(m) \} $ of $ \bN $ 
and some $ \Lambda_* \in \cP(V^{m_\mu}_\times) $. 
By Remark \ref{rem: mupn to eta}, we have 
\begin{align}
\Lambda_0 = \mu^{pn(m)} * \Lambda_{-pn(m)} 
\to \eta * \Lambda_* , 
\label{}
\end{align}
which shows 
\begin{align}
\Lambda_0 = \eta * \Lambda_* = \eta^L * \omega_H * \eta^R * \Lambda_* . 
\label{}
\end{align}
By the same argument as that after \eqref{eq: mu*Lambda = Lambda}, 
we obtain $ \cS(\Lambda_*) \subset W_\mu $. 

Since $ \cS(\eta^R * \Lambda_*) \subset RW_\mu = eW_\mu = GW = CHW $, we have 
\begin{align}
\eta^L * \omega_H * \eta^R * \Lambda_* 
=& \sum_{r=0}^{p-1} \sum_{h \in H} \sum_{\vw \in W} 
(\eta^L * \omega_H * \delta_{\gamma^rh\vw}) \, (\eta^R * \Lambda_*)\{ \gamma^rh\vw \} 
\label{} \\
=& \sum_{r=0}^{p-1} \sum_{\vw \in W} 
(\eta^L * \delta_{\gamma^r} * \omega_H * \delta_{\vw}) \, (\eta^R * \Lambda_*)(\gamma^rH\vw) 
\label{} \\
=& \sum_{r=0}^{p-1} \alpha _r \eta^L * \delta_{\gamma^r} * \omega_H * \Lambda_W^r , 
\label{}
\end{align}
where we set 
\begin{align}
\alpha _r = (\eta^R * \Lambda_*)(\gamma^rHW) 
, \quad 
\Lambda_W^r = \frac{1}{\alpha _r} 
\sum_{\vw \in W} (\eta^R * \Lambda_*)(\gamma^rH\vw) \delta_{\vw} . 
\label{}
\end{align}
We then obtain \eqref{eq: Lambdak} for $ k=0 $. 
By \eqref{eq: mu*...}, 
we obtain \eqref{eq: Lambdak} also for $ k \ge 1 $. 

Let us prove \eqref{eq: Lambdak} for $ k \le -1 $ by induction. 
Suppose \eqref{eq: Lambdak} for a fixed $ k \le 0 $ hold true. 
We want to prove \eqref{eq: Lambdak} for $ k-1 $. 
By the same argument as for $ \Lambda_0 $, 
we have 
\begin{align}
\Lambda_{k-1} = \sum_{r=0}^{p-1} \tilde{\alpha }_r \, \eta^L * \delta_{\gamma^{k-1+r}} * \omega_H * \tilde{\Lambda}^r_W 
\label{}
\end{align}
for some $ \tilde{\Lambda}^0_W,\ldots,\tilde{\Lambda}^{p-1}_W \in \cP(W) $ 
and some constants $ \tilde{\alpha }_0,\ldots,\tilde{\alpha }_{p-1} \ge 0 $ 
such that $ \tilde{\alpha }_0 + \cdots + \tilde{\alpha }_{p-1} = 1 $. 
We then have 
\begin{align}
\Lambda_k 
= \mu * \Lambda_{k-1} 
= \sum_{r=0}^{p-1} \tilde{\alpha }_r \, \eta^L * \delta_{\gamma^{k+r}} * \omega_H * \tilde{\Lambda}^r_W . 
\label{}
\end{align}
Comparing this identity with \eqref{eq: Lambdak} 
and using Proposition \ref{prop: contr W}, we obtain 
$ \alpha_r = \tilde{\alpha }_r $ and $ \Lambda^r_W = \tilde{\Lambda}^r_W $ for $ r=0,1,\ldots,p-1 $. 
We thus obtain \eqref{eq: Lambdak} for $ k-1 $. 
We have proved \eqref{eq: Lambdak} for $ k \le -1 $ by induction. 
}

We now deal with the non-stationary case 
by reducing it to the stationary case. 

\begin{Thm} \label{thm: non-stat}
Suppose the same assumptions of Proposition \ref{prop: contr W} be satisfied. 
Let $ (\bX,N) $ be an $ m_\mu $-particle $ \mu $-evolution 
such that the sequences $ X^1,\ldots,X^{m_\mu} $ are distinct a.s., 
i.e., [$ X^i \neq X^j $ whenever $ i \neq j $] a.s. 
Set $ \Lambda_k(\cdot) := \bP(\bX_k \in \cdot) $ for $ k \in \bZ $. 
Then the following hold: 
\begin{enumerate}

\item[$ \bullet $] 
For any $ k \in \bZ $, the states $ X^1_k,\ldots,X^{m_\mu}_k $ 
are distinct and form an F-clique, a.s. 

\item[$ \bullet $] 
$ (\Lambda_k)_{k \in \bZ} $ satisfies 
the equivalent conditions of Proposition \ref{prop: char Lambda k}. 

\item 
For any fixed $ k \in \bZ $, 
it holds that $ \bX_k \in W_\mu = LGW $ a.s. 
and $ \bX_k^L \dist \eta^L $. 

\item 
$ \bX_k^G = (\gamma^k Y_C)^C U_k $ a.s. for $ k \in \bZ $ 
for some $ C $-valued random variable $ Y_C $ and some $ H $-valued random variables $ U_k $ 
such that $ U_k $ is uniform on $ H $. 

\item 
$ \bX_k^W = \bZ_W $ a.s. for $ k \in \bZ $ 
for some $ W $-valued random variable $ \bZ_W $. 

\item 
If we write $ M^G_j := \bX^G_j (\bX^G_{j-1})^{-1} $ for $ j \in \bZ $ 
and $ M^G_{k,j} := \bX^G_k (\bX^G_j)^{-1} = M^G_k M^G_{k-1} \cdots M^G_{j+1} $ for $ j \le k $, 
we have the following factorization: 
\begin{align}
\bX_j = \bX_j^L (M^G_{k,j})^{-1} (\gamma^k Y_C)^C U_k \bZ_W 
\quad \text{a.s. for $ j \le k $}. 
\label{eq: bXj factorization2}
\end{align}

\item 
A resolution of the observation holds in the sense that 
\begin{align}
\cF^{\bX}_k = \cG^N_k \vee \cF^{\bX}_{-\infty } \vee \sigma(U_k) 
\quad \text{a.s.}, 
\label{eq: resol2}
\end{align}
where 
\begin{align}
\cG^N_k = \sigma \rbra{ \bX^L_j , \ M^G_j : j \le k } \subset \cF^N_k (\subset \sigma(N)) 
\quad \text{a.s.} , 
\label{eq: cGNk cFNk2}
\end{align}
\begin{align}
\text{the three $ \sigma $-fields $ \sigma(N) $, 
$ \cF^{\bX}_{-\infty } $ and $ \sigma(U_k) $ are independent} 
\label{eq: main indep2}
\end{align}
and 
\begin{align}
\cF^{\bX}_{-\infty } 
= \sigma(Y_C , \ \bZ_W) 
\quad \text{a.s.} 
\label{eq: bX-infty2}
\end{align}

\item 
Let $ \alpha _0,\ldots,\alpha _{p-1} $ 
and $ \Lambda_W^0,\ldots,\Lambda_W^{p-1} $ be as in Proposition \ref{prop: char Lambda k}.
Then the joint distribution of $ Y_C $ and $ \bZ_W $ is given as 
\begin{align}
\bP(Y_C=\gamma^r, \ \bZ_W = \vw) 
= \alpha_r \Lambda^r_W \{ \vw \} 
\quad \text{for $ r=0,\ldots,p-1 $ and $ \vw \in W $.} 
\label{}
\end{align}

\end{enumerate}
\end{Thm}

\Proof{
By Proposition \ref{prop: num of distinct} 
and by the assumption that the sequences $ X^1,\ldots,X^{m_\mu} $ are distinct a.s., 
we see that, for every $ k \in \bZ $, 
the states $ X^1_k,\ldots,X^{m_\mu}_k $ are distinct and form an F-clique, a.s. 

We now have $ \bX_k \in V^{m_\mu}_\times $ a.s., 
which shows $ \Lambda_k \in \cP(V^{m_\mu}_\times) $. 
By definition of $ \mu $-evolution, 
we see that Condition (i) of Proposition \ref{prop: char Lambda k} is satisfied. 
Hence we have a representation \eqref{eq: Lambdak}. 

We write $ \omega_W $ for the uniform probability on $ W $ 
and write $ \tilde{\Lambda} = \eta^L * \omega_G * \omega_W $, 
which is a $ \mu $-invariant probability supported on $ W_\mu $. 
Let $ (\tilde{\bX},\tilde{N}) $ under $ \tilde{\bP} $ 
be a stationary $ m_\mu $-particle $ \mu $-evolution 
such that $ \tilde{\bX} $ has $ \tilde{\Lambda} $ as its common law. 
By (vi) of Theorem \ref{thm: main}, we know that 
\begin{align}
\tilde{\bP}(\tilde{Y}_C = \gamma^r , \ \tilde{\bZ}_W = \vw) 
= \frac{1}{p} \cdot \frac{1}{\#(W)} > 0 
\label{}
\end{align}
and so the conditional probability 
\begin{align}
\tilde{\bP}_{\gamma^r,{\vw}}(\cdot) 
:= \tilde{\bP}(\cdot \mid \tilde{Y}_C = \gamma^r , \ \tilde{\bZ}_W = \vw) 
\label{}
\end{align}
is well-defined. 
We then see that $ (\tilde{\bX},\tilde{N}) $ under $ \tilde{\bP}_{\gamma^r,{\vw}} $ 
is a (non-stationary) $ m_\mu $-particle $ \mu $-evolution; 
Since the random variables $ \tilde{Y}_C $ and $ \tilde{\bZ}_W $ are 
$ \cF^{\tilde{\bX}}_{-\infty } $-measurable, 
conditioning by them preserves Markov property \eqref{eq: MP}. 
Note that, for each $ k \in \bZ $, 
the law of $ \tilde{\bX}_k $ under $ \tilde{\bP}_{\gamma^r,{\vw}} $ 
equals to $ \eta^L * \delta_{(\gamma^{k+r})^C} * \omega_H * \delta_{\vw} $. 
Moreover, by \eqref{eq: bXj factorization}, we obtain the following factorization: 
\begin{align}
\tilde{\bX}_j = \tilde{\bX}_j^L (\tilde{M}^G_{k,j})^{-1} (\gamma^{k+r})^C \tilde{U}_k \vw 
\quad \text{$ \tilde{\bP}_{\gamma^r,{\vw}} $-a.s. for $ j \le k $} , 
\label{eq: factrz tildebPiw}
\end{align}
where $ \tilde{M}^G_{k,j} $ and $ \tilde{U}_k $ 
are defined in the same way as in Theorem \ref{thm: main}. 
We then see that Claims (i)-(iv) are satisfied 
for $ (\tilde{\bX},\tilde{N}) $ under $ \tilde{\bP}_{\gamma^r,{\vw}} $. 

Let us check that (v) is satisfied 
for $ (\tilde{\bX},\tilde{N}) $ under $ \tilde{\bP}_{\gamma^r,{\vw}} $. 
By \eqref{eq: cGNk cFNk} of Theorem \ref{thm: main}, 
there exist measurable maps $ \phi_k $ and $ \psi_k $ such that 
$ \tilde{\bX}^L_k = \phi_k(\tilde{N}_j:j \le k) $ 
and $ \tilde{M}^G_k = \psi_k(\tilde{N}_j:j \le k) $, 
$ \tilde{\bP} $-a.s., 
which yields that 
these identities also hold $ \tilde{\bP}_{\gamma^r,{\vw}} $-a.s., 
which shows that Claim \eqref{eq: cGNk cFNk2} is satisfied. 
Since the $ \sigma $-field 
$ \cF^{\tilde{\bX}}_{-\infty } = \sigma(\tilde{Y}_C , \ \tilde{\bZ}_W) $ 
is trivial $ \tilde{\bP}_{\gamma^r,{\vw}} $-a.s., 
we can deduce from the factorization (iv) 
that Claims \eqref{eq: resol2} and \eqref{eq: main indep2} are satisfied.

Define 
\begin{align}
\tilde{\bQ} 
= \sum_{r=0}^{p-1} \alpha _r \sum_{{\vw} \in W} \Lambda^r_W \{ \vw \} 
\, \tilde{\bP}_{\gamma^r,{\vw}} . 
\label{}
\end{align}
We then see that the joint law of $ (\bX,N) $ under $ \bP $ 
equals to that of $ (\tilde{\bX},\tilde{N}) $ under $ \tilde{\bQ} $; 
in fact, they are $ \mu $-evolutions and 
\begin{align}
\tilde{\bQ}(\tilde{\bX}_k \in \cdot) 
=& \sum_{r=0}^{p-1} \alpha _r \sum_{{\vw} \in W} \Lambda^r_W \{ \vw \} 
\, (\eta^L * \delta_{(\gamma^{k+r})^C} * \omega_H * \delta_{\vw}) 
\label{} \\
=& \sum_{r=0}^{p-1} \alpha _r \, \eta^L * \delta_{(\gamma^{k+r})^C} * \omega_H * \Lambda^r_W 
= \Lambda_k 
= \bP(\bX_k \in \cdot) . 
\label{}
\end{align}
We now obtain Claim (vi) 
for $ (\bX,N) $ under $ \bP $. 
We thus derive from \eqref{eq: factrz tildebPiw} the following factorization: 
\begin{align}
\tilde{\bX}_j = \tilde{\bX}_j^L (\tilde{M}^G_{k,j})^{-1} (\gamma^k \tilde{Y}_C)^C \tilde{U}_k \tilde{\bZ}_W 
\quad \text{$ \tilde{\bQ} $-a.s. for $ j \le k $} , 
\label{}
\end{align}
where $ \tilde{Y}_C $ and $ \tilde{\bZ}_W $ 
are defined in the same way as in Theorem \ref{thm: main}. 
We therefore obtain Claims (i)-(v) 
for $ (\bX,N) $ under $ \bP $. 
}

\section{Appendix: The semigroup consisting of mappings} \label{sec: app1}

\Proof[Proof of Proposition \ref{prop: mapping semigroup}]{ 
(i) 
Let us define $ K $ by \eqref{eq: def of K} 
and prove that $ K $ is a minimal ideal of $ S $. 
For $ f \in K $ and $ g,h \in S $, we have 
$ m_S \le \#(gfhV) \le \#(fV) = m_S $, which shows that $ K $ is an ideal of $ S $. 
Let $ I \subset K $ be an ideal of $ S $. 
Let $ f \in K $ and $ g \in I $. Since $ fgf|_{fV} $ is a permutation of $ fV $, 
there exists an integer $ q \ge 1 $ such that $ (fgf)^q $ is identity on $ fV $, 
which implies $ (fgf)^q f = f $. 
Hence 
\begin{align}
f = (fgf)^q f = fg(f(fgf)^{q-1}f) \in SIS \subset I , 
\label{}
\end{align}
which shows $ I=K $, 
and thus $ K $ is the kernel of $ S $. 

(ii) This is obvious. 

(iii) Let $ e $ be a primitive idempotent of $ S $. 
Since $ K $ is completely simple by Proposition \ref{prop: kernel}, 
we may take $ f \in E(K) $. Then $ efe \in SKS \subset K $. 
Since $ efe|_{efeV} $ is a permutation of $ efeV $, 
there exists an integer $ q \ge 1 $ such that $ (efe)^q $ is identity on $ efeV $, 
which yields $ (efe)^{q+1} = efe $. 
If we write $ g := (efe)^{2q} $, we obtain $ eg=ge=g \in E(K) $, 
which implies $ g=e $ by primitivity. Thus we obtain $ e \in E(K) $. 
The converse is obvious since all idempotents of $ K $ are primitive. 

(iv) Let $ f \in Se = LG $ and take $ (x,g) \in L \times G $ 
such that $ f = xg $. 
Since $ g^{-1}f = e $ and $ fe = f $, we have 
[$ fv=fw $ $ \iff $ $ ev=ew $] for all $ v,w \in V $, 
which shows $ \pi(f)=\pi(e) $. 

Conversely, let $ f \in S $ be such that $ \pi(f)=\pi(e) $. 
Then $ \#(fV) = \#(eV) = m_S $, so that $ f \in K $. 
Let $ f = xgy $ with $ (x,g,y) \in L \times G \times R $. 
Since $ \pi(y) = \pi(gy) = \pi(ef) = \pi(f) = \pi(e) $ and $ ye = e $, 
we obtain $ y = e $, so $ f \in Se $. 

(v) Let $ f \in eS = GR $ and take $ (g,y) \in G \times R $ 
such that $ f = gy $. 
Then $ fV = efV \subset eV $. 
Since $ \#(fV) = \#(eV) = m_S $, we have $ fV=eV $. 

Conversely, let $ f \in S $ be such that $ fV=eV $. 
Then $ \#(fV) = \#(eV) = m_S $, so that $ f \in K $. 
Take $ (x,g,y) \in L \times G \times R $ such that $ f = xgy $. 
Note that $ fe = xgye = xg $ and $ x = feg^{-1} = fg^{-1} $. 
Since $ xV = fg^{-1}V \subset fV = eV $, we have $ xV = eV $. 
On one hand, since $ e $ is identity on $ xV = eV $, we have $ exv = xv $ for $ v \in V $. 
On the other hand, since $ ex = e $, we have $ exv = ev $ for $ v \in V $. 
We now obtain $ x=e $, so $ f \in eS $. 

(vi) This is immediate from (iv) and (v), 
since $ G = Se \cap eS $. 
}

\section{Appendix: Another example} \label{sec: app2}

Let 
\begin{align}
V = \cbra{ \pmat{a_1 \\ a_2 \\ a_3} : a_1,a_2,a_3 \in \{ -1,1 \} } . 
\label{}
\end{align}
Let $ D = \{ (1,0),(-1,0),(0,1),(0,-1) \} $ and set 
\begin{align}
S = \cbra{ \pmat{ a_{11} & a_{12} & 0 \\ 
a_{21} & a_{22} & 0 \\ 0 & 0 & b } : 
(a_{11},a_{12}) \in D , \ (a_{21},a_{22}) \in D , \ b \in \{ -1,1 \} }. 
\label{}
\end{align}
Note that $ S $ is a finite semigroup with respect to the usual matrix product. 
In fact, 
\begin{align}
\pmat{ 1 & 0 } \pmat{ \pm 1 & 0 \\ a_{21} & a_{22} } 
= \pmat{ \pm 1 & 0 } 
, \quad 
\pmat{ 1 & 0 } \pmat{ 0 & \pm 1 \\ a_{21} & a_{22} } 
= \pmat{ 0 & \pm 1 } 
\label{}
\end{align}
for $ (a_{21},a_{22}) \in D $, etc. 
We regard an element of $ S $ as a map of $ V $ into itself 
with respect to the usual matrix product. 
Set 
\begin{align}
s_0 = \pmat{ 1 & 0 & 0 \\ 0 & 1 & 0 \\ 0 & 0 & -1 } 
, \quad 
s_1 = \pmat{ -1 & 0 & 0 \\ 0 & 1 & 0 \\ 0 & 0 & -1 } 
, \quad 
s_2 = \pmat{ 0 & 1 & 0 \\ 1 & 0 & 0 \\ 0 & 0 & -1 } 
, \quad 
\label{}
\end{align}
and 
\begin{align}
g = \pmat{ 0 & 1 & 0 \\ 0 & 1 & 0 \\ 0 & 0 & -1 } , 
\label{eq: def of g}
\end{align}
so that $ s_0,s_1,s_2,g \in S $. 
Let $ \mu = (\delta_{s_0} + \delta_{s_1} + \delta_{s_2} + \delta_{g})/4 $ 
be the uniform law on $ \cS(\mu) = \{ s_0,s_1,s_2,g \} $. 
It is easy to see that $ \cS(\mu) $ generates $ S $, i.e., 
$ S = \bigcup_{n=1}^{\infty } \{ s_0,s_1,s_2,g \}^n $. 
Let us apply 
Propositions \ref{prop: IC}, \ref{prop: mapping semigroup} and \ref{prop: contr W}. 

If we write 
\begin{align}
A = \cbra{ \pmat{a_1 \\ a_2 \\ a_3} \in V : a_1=a_2 } 
, \quad 
B = \cbra{ \pmat{a_1 \\ a_2 \\ a_3} \in V : a_1=-a_2 } , 
\label{}
\end{align}
then 
\begin{align}
\begin{cases}
s_0A=A \\
s_0B=B 
\end{cases}
\quad 
\begin{cases}
s_1A=B \\
s_1B=A 
\end{cases}
\quad 
\begin{cases}
s_2A=A \\
s_2B=B 
\end{cases}
\label{}
\end{align}
and $ gV=gA=gB=A $, 
and hence we see that $ m_\mu = m_S = \#(A) = \#(B) = 4 $. 

Set 
\begin{align}
S_+ = \cbra{ \pmat{ a_{11} & a_{12} & 0 \\ 
a_{21} & a_{22} & 0 \\ 0 & 0 & b } \in S : b=1 } 
, \quad 
S_- = \cbra{ \pmat{ a_{11} & a_{12} & 0 \\ 
a_{21} & a_{22} & 0 \\ 0 & 0 & b } \in S : b=-1 } . 
\label{eq: def of S-}
\end{align}
Since $ \cS(\mu) \subset S_- $, we have 
$ \cS(\mu^2) = \cS(\mu) \cS(\mu) \subset S_- S_- = S_+ $. 
Since $ \cS(\mu^2) $ generates $ S_+ $ and contains the identity map, 
we see that the left or right random walk on $ S_+ $ 
whose steps have law $ \mu^2 $ is aperiodic, 
whereas the random walk on $ S $ whose steps have law $ \mu $ is not. 
Hence we obtain $ p=2 $, and consequently 
the sequence $ \{ \mu^{2n} \}_{n=1}^{\infty } $ converges to $ \eta $. 

Let $ K = \cS(\nu) $ and $ K_+ = \cS(\eta) $ 
denote the kernels of $ S $ and $ S_+ $, respectively. 
Then 
\begin{align}
K =& \{ s \in S : \#(sV) = 4 \} 
= \cbra{ \pmat{ a_{11} & 0 & 0 \\ a_{21} & 0 & 0 \\ 0 & 0 & b }, 
\pmat{ 0 & a_{11} & 0 \\ 0 & a_{21} & 0 \\ 0 & 0 & b } \in S } 
, \label{} \\
K_+ =& \{ s \in S_+ : \#(sV) = 4 \} 
= \cbra{ \pmat{ a_{11} & 0 & 0 \\ a_{21} & 0 & 0 \\ 0 & 0 & 1 }, 
\pmat{ 0 & a_{11} & 0 \\ 0 & a_{21} & 0 \\ 0 & 0 & 1 } \in S } . 
\label{}
\end{align}
Set 
\begin{align}
e := g^2 = \pmat{ 0 & 1 & 0 \\ 0 & 1 & 0 \\ 0 & 0 & 1 } \in S_+ , 
\label{}
\end{align}
which is an idempotent of $ S_+ $. 
Since $ \#(eV) = 4 $, we have $ e \in K_+ $. 
Let us determine the Rees decompositions $ K=LGR $ and $ K_+=LHR $ at $ e \in E(K_+) $. 
Set 
\begin{align}
f = \pmat{ 0 & -1 & 0 \\ 0 & 1 & 0 \\ 0 & 0 & 1 } 
, \quad 
k = \pmat{ 1 & 0 & 0 \\ 1 & 0 & 0 \\ 0 & 0 & 1 } . 
\label{}
\end{align}
Then we have 
\begin{align}
L =& \{ s \in E(K_+) : \pi(s) = \pi(e) \} = \{ e,f \} , 
\label{} \\
R =& \{ s \in E(K_+) : sV = eV \} = \{ e,k \} , 
\label{} \\
H =& \{ s \in K_+ : \pi(s) = \pi(e) , \ sV=eV \} = \{ e,-g \} , 
\label{} \\
G =& \{ s \in K : \pi(s) = \pi(e) , \ sV=eV \} = \{ e,-e,g,-g \} . 
\label{}
\end{align}
We now see that we may choose $ \gamma = g $, so that $ C = \{ e,g \} $ and $ G = CH $. 
We have the following multiplication tables 
(the table of $ ab $ for $ a $ and $ b $): 
\begin{align}
\text{\begin{tabular}{c|cccc}
$ a \backslash b $ & $ e $ & $ f $ & $ g $ & $ k $ 
\\ \hline
$ e $ & $ e $ & $ e $ & $ g $ & $ k $ 
\\
$ f $ & $ f $ & $ f $ & $ fg $ & $ fk $ 
\\
$ g $ & $ g $ & $ g $ & $ e $ & $ gk $ 
\\
$ k $ & $ e $ & $ -g $ & $ g $ & $ k $ 
\end{tabular}}
\qquad 
\text{\begin{tabular}{c|cc}
$ a \backslash b $ & $ e $ & $ f $ 
\\ \hline
$ s_0 $ & $ g $ & $ fg $ 
\\
$ s_1 $ & $ fg $ & $ g $ 
\\
$ s_2 $ & $ g $ & $ -f $ 
\\
$ g $ & $ g $ & $ g $ 
\end{tabular}}
\qquad 
\text{\begin{tabular}{c|cccc}
$ a \backslash b $ & $ s_0 $ & $ s_1 $ & $ s_2 $ & $ g $ 
\\ \hline
$ e $ & $ g $ & $ g $ & $ gk $ & $ g $ 
\\
$ k $ & $ gk $ & $ -k $ & $ g $ & $ g $ 
\end{tabular}}
\label{eq: mtable}
\end{align}
Note that $ -f = fg(-g) $ and $ -k = g(-g)k $. 

Let us compute $ \eta^L $. Let $ \eta^L = \alpha \delta_e + \beta \delta_f $ 
for some $ \alpha ,\beta >0 $ with $ \alpha + \beta = 1 $. 
Note that 
\begin{align}
\mu * \eta^L * \omega_H 
=& \mu * \eta^L * \omega_H * \eta^R * \omega_H 
= \mu * \eta * \omega_H 
\label{} \\
=& \eta^L * \delta_g * \omega_H * \eta^R * \omega_H 
= \eta^L * \delta_g * \omega_H . 
\label{}
\end{align}
On the one hand, we have 
\begin{align}
\eta^L * \delta_g * \omega_H 
= \rbra{ \alpha \delta_e + \beta \delta_f } * \delta_g * \omega_H . 
\label{}
\end{align}
On the other hand, 
by \eqref{eq: mtable}, we have 
\begin{align}
\mu * \eta^L 
= \frac{1}{4} \rbra{ \delta_{s_0} + \delta_{s_1} + \delta_{s_2} + \delta_{g} } 
* \rbra{ \alpha \delta_e + \beta \delta_f } 
= \frac{3\alpha + 2\beta }{4} \delta_g + \frac{\alpha + \beta}{4} \delta_{fg} 
+ \frac{\beta}{4} \delta_{(-f)} . 
\label{}
\end{align}
Since $ -f = fg(-g) $ and $ -g \in H $, we have 
$ \delta_{-f} * \omega_H = \delta_{fg} * \omega_H $, so 
\begin{align}
\mu * \eta * \omega_H 
= \mu * \eta^L * \omega_H 
= \rbra{ \frac{3\alpha + 2\beta }{4} \delta_e + \frac{\alpha + 2\beta}{4} \delta_f } 
* \delta_g * \omega_H . 
\label{}
\end{align}
Hence we obtain $ \alpha = 2/3 $ and $ \beta = 1/3 $, so that 
\begin{align}
\eta^L = \frac{2}{3} \delta_e + \frac{1}{3} \delta_f . 
\label{}
\end{align}
By a similar argument, using the identities $ -k = g(-g)k $, 
$ \omega_H * \delta_g = \delta_g * \omega_H $ and 
\begin{align}
\omega_H * \eta^R * \mu 
= \delta_g * \omega_H * \eta^R , 
\label{}
\end{align}
we obtain 
\begin{align}
\eta^R = \frac{2}{3} \delta_e + \frac{1}{3} \delta_k . 
\label{}
\end{align}

Note that $ eV = \{ v_1,v_2,v_3,v_4 \} $ 
and $ fV = \{ fv_1,fv_2,fv_3,fv_4 \} $, 
where 
\begin{align}
v_1 = \pmat{1\\1\\1} 
, \quad 
v_2 = \pmat{-1\\-1\\1} 
, \quad 
v_3 = \pmat{1\\1\\-1} 
, \quad 
v_4 = \pmat{-1\\-1\\-1} . 
\label{}
\end{align}
By \eqref{eq: Snu Fcliques} and \eqref{eq: Wmu rep}, we see that 
\begin{align}
W_\mu = \{ (x^1,x^2,x^3,x^4) : 
\text{a permutation of $ (v_1,v_2,v_3,v_4) $ or $ (fv_1,fv_2,fv_3,fv_4) $} \} 
\label{}
\end{align}
and 
\begin{align}
eW_\mu = \{ (x^1,x^2,x^3,x^4) : 
\text{a permutation of $ (v_1,v_2,v_3,v_4) $} \} . 
\label{}
\end{align}
We have the following multiplication table 
(the table of $ sv $ for $ s $ and $ v $): 
\begin{align}
\text{\begin{tabular}{c|cccc}
$ s \backslash v $ & $ v_1 $ & $ v_2 $ & $ v_3 $ & $ v_4 $ 
\\ \hline
$ e $ & $ v_1 $ & $ v_2 $ & $ v_3 $ & $ v_4 $ 
\\
$ -e $ & $ v_4 $ & $ v_3 $ & $ v_2 $ & $ v_1 $ 
\\
$ g $ & $ v_3 $ & $ v_4 $ & $ v_1 $ & $ v_2 $ 
\\
$ -g $ & $ v_2 $ & $ v_1 $ & $ v_4 $ & $ v_3 $ 
\end{tabular}}
\label{}
\end{align}
From this table, we see that we may take a set $ W $ as 
\begin{align}
W = \{ (x^1,x^2,x^3,x^4) : \text{$ x^4=v_4 $ and $ (x^1,x^2,x^3) $ 
is a permutation of $ (v_1,v_2,v_3) $ } \} , 
\label{}
\end{align}
which is a minimal subset of $ W_\mu $ such that $ eW_\mu = GW $.

\noindent
{\bf Data Availability} 
Data sharing is not applicable to this article as no datasets were generated or analyzed 
during the current study.

\noindent
{\bf Conflict of interest} 
The authors have no competing interests to declare that are relevant to the content 
of this article.

\def\cprime{$'$} \def\cprime{$'$}

\end{document}